\documentclass[12pt]{amsart}

\usepackage{graphicx}
\usepackage{amsfonts}
\usepackage{epsf}
\usepackage{amssymb}
\usepackage{amsmath}
\usepackage{amscd}

\newtheorem{theorem}{Theorem}[section]
\newtheorem{proposition}[theorem]{Proposition}
\newtheorem{lemma}[theorem]{Lemma}
\newtheorem{definition}[theorem]{Definition}

\newtheorem{corollary}[theorem]{Corollary}
 
\newtheorem{remark}[theorem]{Remark} 
\newtheorem{example}[theorem]{Example}

\newcommand{\kerr}{\mbox{Ker} }

\newcommand{\s}{\mathfrak{s}}
\renewcommand{\t}{\mathfrak{t}}
\newcommand{\imm}{\mbox{Im} }

\input xy
\xyoption{all}

\begin{document}
\thanks{The author gratefully acknowledges support from the NSF through grant DMS-0709625.}
\title{Concordance invariants from higher order covers}
\begin{abstract}
{We generalize the Manolescu-Owens smooth concordance invariant $\delta(K)$ of knots $K\subset S^3$ 
to invariants $\delta_{p^n}(K)$ obtained by considering covers of order $p^n$, with $p$ a prime. Our main result shows that for any prime $p\ne 2$, the thus obtained homomorphism $\oplus _{n\in \mathbb{N}} \delta_{p^n}$ from the smooth concordance group to $\mathbb{Z}^\infty$ has infinite rank. We also show that unlike $\delta$, these new invariants  typically are not multiples of the knot signature, even for alternating knots. A significant  portion of the article is devoted to exploring examples.}
\end{abstract}
\author{Stanislav Jabuka}
\address{Department of Mathematics and Statistics, University of Nevada, Reno, NV 89557}
\email{jabuka@unr.edu}
\maketitle
\section{Introduction}
Many of the recent advances in our understanding of the smooth knot concordance group $\mathcal C$
have been driven by the advents of two theories: Heegaard Floer homology  and Khovanov homology. 
The former, discovered and developed by P. Ozsv\'ath and Z.  Szab\'o in a series of 
beautiful papers \cite{peter1, peter2, peter4}, has grown into a comprehensive package of invariants of low dimensional 
manifolds, including invariants of nullhomologous knots in arbitrary 3-manifolds. The latter has been discovered by M.  Khovanov in \cite{khovanov1} and further developed by Khovanov and L. Rosansky in \cite{khovanov2}, and is at present limited in  scope to providing invariants for knots in $S^3$.  

Despite the rather different approaches taken in these two theories, they share a surprising 
amount of formal properties. One such similarity is a pair of spectral sequences associated to a knot $K\subset S^3$.  The $E^2$ terms of these sequences are the knot Floer homology group $\widehat{HFK}(K)$ and 
the Khovanov homology group $Kh(K)$ respectively, while their $E^\infty$ terms are rather standard groups, 
namely $\mathbb{Q}$ and $\mathbb{Q}^2$ (when working with rational coefficients).  These spectral sequences 
have been exploited in \cite{peter11, rasmussen1} and \cite{rasmussen2} to define two epimorphisms $\tau, s : \mathcal C \to \mathbb{Z}$ from the smooth knot concordance group $\mathcal C$ to the integers. The two invariants exhibit a number of similar features and have at first been conjecture to be equal (indeed they agree for all quasi-alternating knots \cite{ozsvath-manolescu}) until the surprising article of M. Hedden and P. Ording \cite{matt-phil} disproved this notion. Both invariants have been found to be powerful obstructions to smooth sliceness. 

In \cite{manolescu-owens} C. Manolescu and B. Owens defined yet another homomorphism $\delta : \mathcal C \to \mathbb{Z}$ by exploiting a different feature of Heegaard Floer homology (thus far unparalleled in Khovanov homology). Namely, the Heegaard Floer homology package associates to a pair $(Y,\s)$ consisting of a 3-manifold $Y$ and a torsion spin$^c$-structure $\s$ on $Y$ (by which we mean that $c_1(\s)$ is torsion in $H^2(Y;\mathbb{Z})$) a rational number $d(Y,\s)$ known as a {\em correction term}, cf. \cite{peter4}. Manolescu and Owens define $\delta (K) = 2 \, d(Y_2(K),\s_0)$ where $Y_2(K)$ is the 2-fold cover of $S^3$ with branching set 
$K$ and $\s_0 \in Spin^c(Y_2(K))$ is the unique spin-structure on $Y_2(K)$. Using properties 
of the correction terms, they show that $\delta$ descends to an epimorphism from $\mathcal C$ to 
$\mathbb{Z}$. 

The utility of the three homomorphisms $\tau, s, \delta:\mathcal C \to \mathbb{Z}$ has been exploited by 
many authors and has led to substantial progress and new results about smooth knot concordance, see for 
example \cite{jabukaglasnik} for a survey. We would like to mention that $\tau, s, \delta$ have been found to be linearly independent
epimorphisms, in fact they remain so even when restricted to the set of topologically slice knots as demonstrated by C. Livingston in \cite{chuck2}. Finally, we point out that for alternating knots, all three of $\tau, s, \delta$ agree with the knot signature, up to a multiplicative constant.  
\vskip3mm
The compelling success of $\tau, s$ and $\delta$ in addressing knot concordance matters is the motivation for the present work. Specifically, the goal of this article is to introduce additional homomorphisms $\delta_{p^n}:\mathcal C \to \mathbb{Z}$ parametrized by 
a pair of positive integers $(p,n)$ of which $p$ is prime. Our construction of $\delta _{p^n}$ exploits the Manolescu-Owens definition of $\delta$ (which in our notation corresponds to $\delta_2$) by using the $p^n$-fold cover of $S^3$ branched over $K$ rather than the $2$-fold cover used in \cite{manolescu-owens}.   
\begin{definition} \label{ofdeltas}
Let $K\subset S^3$ be a knot and $p$ a prime integer. Let $Y_{p^n}(K)$ be the $p^n$-fold 
cover of $S^3$ with branching set $K$. We define the integer $\delta_{p^n}(K)$ as 
$$\delta_{p^n}(K) = 2\,  d(Y_{p^n}(K),\s_0)$$
with $\s_0=\s_0(K,p^n)  \in Spin^c(Y_{p^n})$ being a spin-structure determined by $p^n$ and $K$ (see section \ref{spincstructures-section} for a definition of $\s_0$). We shall refer to $\s_0$ as the canonical spin-structure associated to $(K,p^n)$. 
\end{definition}
It is not immediately clear that the thus defined $\delta_{p^n}(K)$ is indeed an integer since correction terms 
typically take on rational, non-integral values. The normalization factor of $2$ is introduced precisely for that purpose and indeed renders $\delta_{p^n}(K)$ integral. This is demonstrated in section \ref{linear-independence}.  
With this definition in mind, the main result of this article is the next theorem. 
\begin{theorem} \label{main}
For each pair of positive integers $(p,n)$ with $p$ prime, $\delta_{p^n}$ descends to a group homomorphism 
$\delta _{p^n}:\mathcal C \to \mathbb{Z}$ from the smooth knot concordance group $\mathcal C$ to the integers. For a fixed prime $p\ne 2$, the homomorphism  
$$\bigoplus _{n=1}^\infty  \delta _{p^{n}}: \mathcal C \to \mathbb{Z}^\infty $$
is of infinite rank. Moreover $\delta _{2}\oplus \delta _{4}\oplus \delta _{8}\oplus \delta _{16} :\mathcal C \to \mathbb{Z}^4$ is of rank $4$.  
\end{theorem}
The much stronger statement of the above theorem in the case of $p\ne 2$ is based on an understanding of 
when the integers $p^{2^n} + 1$ are prime powers. We give an answer 
to this question in section \ref{linear-independence} for $p\ne 2$. When $p=2$, then numbers $2^{2^n}+1$ are known as {\em Fermat numbers}  and it is an open question in number theory which among them are prime powers. 
It is likely that $\oplus _{n=1}^\infty \delta_{2^n}:\mathcal C \to \mathbb Z^\infty$ 
is also of infinite rank but our proof doesn't apply to this case. No effort was made to maximize 
the rank of $\oplus _{n=1}^\infty \delta_{2^n}$; the rank of $4$ from theorem \ref{main} follows from examples 
\ref{ex1}, \ref{ex5}, \ref{ex6} and \ref{ex7} below. A slightly stronger statement of theorem \ref{main} is given in theorem \ref{mainprime}. 

As already mentioned above, all three of $\tau(K), s(K), \delta(K)$ agree with the signature $\sigma (K)$ (up to a constant scalar factor) for the case of alternating knots $K$. To investigate the relation of the various $\delta_{p^n}$ to the knot signature, we turn to examples.  
%
%
%
\subsection{Examples}
A significant portion of the article is devoted to calculating the invariants $\delta_{p^n}(K)$ for  concrete knots $K$. 
We present 3 different computational techniques (discussed in sections \ref{examples1-section}, \ref{examples2-section} and  \ref{trefoil-section}), for now however we contend ourselves with merely stating the results of these calculations. We let $T_{(a,b)}$ denote the $(a,b)$ torus knot. Throughout we assume that $a,b \in \mathbb{N}$ are relatively prime. 
\begin{example}[The right-handed trefoil $T_{(2,3)}$] \label{ex1}
For those integers $n>0$ for which either of $6n\pm1$ is a prime power, we get 
$$\delta_{6n-1}(T_{(2,3)}) = 4 \quad \quad \mbox{ and } \quad \quad \delta_{6n+1}(T_{(2,3)}) = 0$$
For example $n=1, 2, 3,4,5$ give such prime powers. Additionally
$$\delta_{2^n}(T_{(2,3)}) = \left\{
\begin{array}{cl}
3 & \quad ; \quad \mbox{if $n=2k$ and $k\ge 1$.} \cr
& \cr
1 & \quad ; \quad \mbox{if $n=2k+1$ and $k\ge 0$.}
\end{array}
\right.
$$
\end{example}
\begin{example}[The torus knot $T_{(2,5)}$] \label{ex2} 
For any integer $n>0$ for which $10n\pm 1$ is a prime power (e.g. $n=1,3,6,8,24$ etc.), we get 
$$\delta_{10n-1}(T_{(2,5)}) = 4 \quad \quad \mbox{ and } \quad \quad \delta_{10n+1}(T_{(2,5)}) = 0$$
In addition we also get  
\begin{center}
\begin{tabular}{|c||c|c|c|c|} \hline 
$n$                             & $3$  & $7$  & $13$ & $17$  \cr \hline \hline 
$\delta_n(T_{(2,5)})$  & $4$  & $0$   & $4$ & $0$     \cr \hline 
\end{tabular}
\end{center}
\vskip1mm
\end{example}
\begin{example}[The torus knot $T_{(2,9)}$] \label{ex3}
For any integer $n>0$ making either of $18n\pm 1$ a prime power, gives  
$$\delta_{18n-1}(T_{(2,9)}) = 8 \quad \quad \mbox{ and } \quad \quad \delta_{18n+1}(T_{(2,9)}) = 0$$
Additionally we obtain
\begin{center}
\begin{tabular}{|c||c|c|c|c|} \hline 
$n$                             & $5$  & $7$  & $11$ & $13$    \cr \hline \hline 
$\delta_n(T_{(2,9)})$  & $4$  & $4$  & $0$  & $0$     \cr \hline 
\end{tabular}
\end{center}
\vskip1mm
\end{example}
\begin{example}[The torus knot $T_{(3,4)}$] \label{ex4}
For any integer $n>0$ with either of $12n\pm 1$ a prime power, we find that  
$$\delta_{12n-1}(T_{(3,4)}) = 4 \quad \quad \mbox{ and } \quad \quad \delta_{12n+1}(T_{(3,4)}) = 0$$
Likewise, we find 
\begin{center}
\begin{tabular}{|c||c|c|c|c|} \hline 
$n$                             & $5$  & $7$  & $17$ & $19$    \cr \hline \hline 
$\delta_n(T_{(3,4)})$  & $0$  & $4$  & $0$  & $4$     \cr \hline 
\end{tabular}
\end{center}
\vskip1mm
\end{example}

\begin{example}[The torus knot $T_{(3,5)}$] \label{ex5}
For any integer $n>0$ rendering $15n\pm 1$ a prime power, one finds that 
$$\delta_{15n-1}(T_{(3,5)}) = 8 \quad \quad \mbox{ and } \quad \quad \delta_{15n+1}(T_{(3,5)}) = 0$$
Other examples are

\begin{center}
\begin{tabular}{|c||c|c|c|c|c|c|c|c|} \hline 
$n$                             & $2$  & $4$  & $7$ & $8$& $11$& $13$   & $17$ & $19$  \cr \hline \hline 
$\delta_n(T_{(3,5)})$  & $4$  & $0$  & $4$ & $0$& $4$ & $4$          & $4$ & $0$     \cr \hline 
\end{tabular}
\end{center}
\vskip1mm
\end{example}
\begin{example}[The torus knot $T_{(3,7)}$] \label{ex6}
For any integer $n>0$ rendering $21n\pm 1$ a prime power, one finds that 
$$\delta_{21n-1}(T_{(3,7)}) = 8 \quad \quad \mbox{ and } \quad \quad \delta_{21n+1}(T_{(3,7)}) = 0$$
Similarly we find that 

\begin{center}
\begin{tabular}{|c||c|c|c|c|c|} \hline 
$n$                             & $2$  & $4$  & $5$ & $8$& $16$  \cr \hline \hline 
$\delta_n(T_{(3,7)})$  & $0$  & $4$  & $4$ & $0$& $0$     \cr \hline 
\end{tabular}
\end{center}
\vskip1mm
\end{example}
\begin{example}[The torus knot $T_{(5,7)}$] \label{ex7}
For any integer $n>0$ for which $35n\pm 1$ is a prime power, one finds that 
$$\delta_{35n-1}(T_{(5,7)}) = 16 \quad \quad \mbox{ and } \quad \quad \delta_{35n+1}(T_{(5,7)}) = 0$$
In addition to these we also find 

\begin{center}
\begin{tabular}{|c||c|c|c|c|} \hline 
$n$                             & $2$  & $4$  &  $8$& $16$  \cr \hline \hline 
$\delta_n(T_{(5,7)})$  & $0$  & $0$  & $8$  & $8$       \cr \hline 
\end{tabular}
\end{center}
\vskip1mm
\end{example}
\begin{example}[The torus knot $T_{(5,9)}$] \label{ex8}
For any integer $n>0$ for which $45n\pm 1$ is a prime power, one finds that 
$$\delta_{35n-1}(T_{(5,9)}) =  \quad \quad \mbox{ and } \quad \quad \delta_{35n+1}(T_{(5,9)}) = 0$$
and for $n=2,4,8,16$ we obtain

\begin{center}
\begin{tabular}{|c||c|c|c|c|} \hline 
$n$                             & $2$  & $4$  &  $8$& $16$  \cr \hline \hline 
$\delta_n(T_{(5,9)})$  & $4$  & $4$  & $0$  & $0$       \cr \hline 
\end{tabular}
\end{center}
\vskip1mm
\end{example}
Example \ref{ex1} implies that
\begin{corollary}
For each integer $k\ge 0$, the homomorphism $\delta _{2^{2k+1}}:\mathcal C \to \mathbb{Z}$ is surjective. 
\end{corollary}
\begin{remark} \label{notsign}
As mentioned above, the homomorphisms $\tau, s, \delta : \mathcal C \to \mathbb{Z}$ agree with the signature of the knot (up to a multiplicative constant) for all alternating knots (see \cite{ozsvath-manolescu} for a stronger statement for quasi-alternating knots in the case of $\tau$ and $s$). The above 
examples illustrate that this is not the case for all $\delta _{p^n}$. For example 
$$\delta _7(T_{(2,3)}) = 0 \quad \quad \quad \quad \delta _7(T_{(2,5)}) = 0 \quad \quad \quad \quad \delta _7(T_{(2,9)}) = 4 $$ 
while the signatures of these knots are 
$$ \sigma (T_{(2,3)}) = -2  \quad \quad \quad \quad  \sigma (T_{(2,5)}) = -4 \quad \quad \quad \quad  \sigma (T_{(2,9)}) = - 8$$
%
\end{remark}
\subsection{Organization}
The remainder of the article is organized as follows. Sections \ref{spincstructures-section} and \ref{backgroundmaterial}  provide background material on spin$^c$-structures (with an emphasis on branched covers) and Heegaard Floer homology. Definition \ref{ofs0} in section \ref{spincstructures-section} specifies the canonical spin-structure $\s_0(K,p^n)$ alluded to in definition \ref{ofdeltas}.  Section \ref{linear-independence} is devoted to the proof of theorem \ref{main} and its slightly strengthened version, theorem \ref{mainprime}. Sections \ref{examples1-section} -- \ref{trefoil-section} are devoted 
to exploring computational techniques. In particular, the results from examples \ref{ex1} -- \ref{ex8} follow directly 
from the discussions in those sections. 
\vskip3mm
\noindent {\bf Acknowledgement } I would like to thank Ron Fintushel for a helpful email exchange. 
\section{Spin$^c$-structures} \label{spincstructures-section}
%
%
%
\subsection{Spin$^c$-structures on three and four manifolds}  \label{spinc34}
This section discusses spin$^c$-structures on $3$ and $4$ manifolds. Our exposition 
largely follows that from chapter 11 in Turaev's book \cite{turaev2}, see also \cite{turaev1}. To begin with, 
recall that the groups $Spin(n)$ for $n\ge 3$ are defined to be the universal covering spaces of $SO(n)$
(it is well known that $\pi_1(SO(n))\cong \mathbb{Z}_2$ for $n\ge 3$, see for example \cite{hatcher}) and 
$Spin^c(n)$ is defined as    
$$ Spin^c(n) = (Spin(n) \times U(1))/\mathbb{Z}_2$$
where $\mathbb{Z}_2=\{ \pm 1\}$ acts by diagonal multiplication. In the cases of $n=3, 4$ one obtains group isomorphisms $Spin^c(3) \cong U(2)$ and $Spin^c(4) \cong U(2) \times U(2)$. Similarly, there are 
diffeomorphisms $SO(3) \cong \mathbb{RP}^3$ and $SO(4) \cong S^3 \times \mathbb{RP}^3$, cf. \cite{hatcher}. For later use we point out that 
\begin{equation} \nonumber
H^2(SO(n);\mathbb{Z}) \cong \mathbb{Z}_2 \quad  \quad \mbox{ and } \quad \quad H^1(SO(n);\mathbb{Z}_2) \cong \mathbb{Z}_2 \quad \quad \mbox{ for } n=3,4.
\end{equation}
If one thinks of $Spin^c(n)\to SO(n)$ as a $U(1)$-bundle, then in both the cases of $n=3,4$, the first Chern 
class of the bundle corresponds to the nontrivial element of $H^2(SO(n);\mathbb{Z})$.

In the following we let $Y$ be a 3-manifold and $X$ a 4-manifold, both possibly with boundary. 
All manifolds are always assumed to be smooth, compact and oriented. For convenience we endow our
manifolds with a Riemannian metric which we assume to be a product metric in a collar neighborhood of the 
boundary. By {\em the frame bundle} we shall mean the bundle of oriented orthonormal frames and we 
shall denote it by $Fr_Y$ or $Fr_X$. These are, of course, principal $SO(3)$ and $SO(4)$ bundles respectively. 

A spin$^c$-structure $\s$ on $Y$ is a principal $Spin^c(3)$ bundle $P_{Spin^c(3)}\to Y$ together with an bundle map $\alpha :P_{Spin^c(3)}\to Fr_Y$  which fiberwise restricts to give the above map $Spin^c(3) \to SO(3)$
(i.e. projection from $Spin^c(3)$ to $Spin(3)$, followed by the covering map to $SO(3)$). An alternative and 
equivalent point of view is to think of $\s$ as a $U(1)$ bundle over $Fr_Y$ which restricts over each fiber 
$SO(3)\hookrightarrow Fr_Y$ to give the unique nontrivial $U(1)$-bundle over $SO(3)$. Said differently, we can define a spin$^c$-structure $\s$ on $Y$ as an element from $H^2(Fr_Y;\mathbb{Z})$ which on each fiber $SO(3)\hookrightarrow Fr_Y$ restricts to the unique nontrivial element of $H^2(SO(3);\mathbb{Z})\cong \mathbb{Z}_2$. 
Fixing an orthogonal trivialization of $TY$, we obtain the diffeomorphism $Fr_Y\cong Y\times SO(3)$ and therefore 
$$H^2(Fr_Y;\mathbb{Z})\cong H^2(Y;\mathbb{Z})\oplus H^2(SO(3);\mathbb{Z})\cong  H^2(Y;\mathbb{Z})\oplus 
\mathbb{Z}_2$$
A spin$^c$-structure $\s$ on $Y$ is thus an element of $H^2(Fr_Y;\mathbb{Z})$ with nontrivial second 
coordinate in this decomposition. There is an obvious action of $H^2(Y;\mathbb{Z})$ on $Spin^c(3)$ given by the pullback map on second cohomology induced by the bundle map $Fr_Y\to Y$ along with addition in $H^2(Fr_Y;\mathbb{Z})$. This action is obviously free and transitive revealing that $Spin^c(Y)$ is an $H^2(Y;\mathbb{Z})$ affine space. 

$Spin^c(Y)$ comes equipped with an involution sending and element in $H^2(Fr_Y;\mathbb{Z})$ to its negative. 
We shall denote this map by $\s \mapsto \bar \s $ and call $\bar \s$ the {\em conjugate spin$^c$-structure of $\s$}.
Finally, note that if $Y'\subset Y$ is a codimension zero submanifold, there is a natural restriction induced map $Spin^c(Y) \to Spin^c(Y')$. In applications below, $Y'$ will typically be the complement of a tubular neighborhood of a knot in $Y$. 

The same arguments apply verbatim to the $4$-manifold $X$ as well (though the case of $X$ closed and $TX$ nontrivial requires a slightly different argument to show that $H^2(Fr_X;\mathbb{Z})\cong H^2(X;\mathbb{Z})\oplus H^2(SO(4);\mathbb{Z})$, it involves a choice of a trivialization of $TX$ over $X$ minus a $4$-ball).   
\vskip2mm
A spin-structure on $Y$ is a principal $Spin(3)$-bundle $P_{Spin(3)}\to Y$ with a bundle map to 
$Fr_Y$ which fiberwise restricts to the double covering map $Spin(3) \to SO(3)$. We denote the set of 
spin-structures on $Y$ by $Spin(Y)$. Thinking of $P_{Spin(3)}$ as 
a (real) line bundle over $Fr_Y$, we can alternatively define a spin-structure on $Y$ as an element of 
$H^1(Fr_Y;\mathbb{Z}_2)$ which on each fiber $SO(3)$ restricts to the nontrivial element of
$H^1(SO(3);\mathbb{Z}_2)$.  Just as with the case of spin$^c$-structures, the obvious action of $H^1(Y;\mathbb{Z}_2)$ on $Spin(Y)$ gives the latter the structure of an affine $H^1(Y;\mathbb{Z}_2)$-space. 
Identical definitions and properties apply to $Spin(X)$, the space of spin-structures on $X$. 

There is a natural homomorphism of affine spaces $Spin(Y) \to Spin^c(Y)$ given by the Bockstein 
map $H^1(Fr_Y;\mathbb{Z}_2)\to H^2(Fr_Y;\mathbb{Z})$ associated to the exact sequence 
$0\to \mathbb{Z} \stackrel{\cdot 2}{\to} \mathbb{Z} \to \mathbb{Z}_2\to 0$. Thus, under a compatible choice 
of origins in $Spin(Y)$ and $Spin^c(Y)$, a spin$^c$-structure $\s \in H^2(Y;\mathbb{Z})$ is a spin-structure 
if and only if 
\begin{equation} \nonumber 
\s \in \imm \left( H^1(Y;\mathbb{Z}_2)\to H^2(Y;\mathbb{Z})  \right) = \kerr \left( H^2(Y;\mathbb{Z}) \stackrel{\cdot 2}{\to} H^2(Y;\mathbb{Z}) \right) 
\end{equation}
Note that distinct spin-structures may descent to give the same spin$^c$-structure. 
Similar considerations apply to $X$. 
\subsection{Spin$^c$-structures on connected sums} 
Let $Y_0$ and $Y_1$ be two closed 3-manifolds and let $Y=Y_0\#Y_1$ be their connected sum. 
Let $B_i\subset Y_i$ be the 3-balls used to perform the connected sum. 
An easy exercise in homological algebra reveals that the restriction maps $H^2(Y_i;\mathbb{Z})\to 
H^2(Y_i-B_i;\mathbb{Z})$ are isomorphisms giving rise to the isomorphism (of affine spaces) 
$Spin^c(Y_i) \to Spin^c(Y_i-B_i)$. We shall utilize these isomorphisms to identify spin$^c$-structures on $Y_i$ and 
$Y_i-B_i$.  

Yet another easy exercise shows that the restriction map 
$$H^2(Y;\mathbb{Z}) \to H^2(Y_1-B_1;\mathbb{Z}) \oplus H^2(Y_2-B_2;\mathbb{Z})$$ 
is an isomorphism as well allowing us to identify $Spin^c(Y)$ with $Spin^c(Y_1-B_1)\times Spin^c(Y_2-B_2)$.  
Putting these two observations together yields the isomorphism 
\begin{equation} \label{addingspinc}
Spin^c(Y_1)\times Spin^c(Y_2) \cong Spin^c(Y_1\#Y_2) \quad \quad \quad \quad (\s_1, \s_2) \mapsto \s_1 \# \s_2
\end{equation}
%
\subsection{Spin$^c$-structures on branched covers}  \
We now turn our attention to spin$^c$-structures on branched covers. 
We adopt the convention that whenever we consider a smooth map $f:Y_0\to Y_1$ which is a local diffeomorphism, the Riemannian metric on $Y_0$ shall be the one obtained from the metric on $Y_1$ via pullback through $f$. 
 
By way of notation, let $K\subset S^3$ be a knot and let $Y_{p^n}(K)$ be the $p^n$-fold branched cover of $S^3$ with branching set $K$. We always assume that $p$ is prime so that $Y_{p^n}(K)$ is a 
rational homology sphere (for a nice proof of this fact see \cite{chuck1}). We let $f:Y_{p^n}(K) \to S^3$ denote the branch covering map and we let $K'=f^{-1}(K)$. We shall write $N(K)$ and $N(K')$ to denote tubular neighborhoods of $K$ and $K'$ respectively. 

Since $f:(Y_{p^n}(K)-N(K')) \to (S^3-N(K))$ is a local diffeomorphism, it induces a push-forward map 
$f_*:Fr_{Y_{p^n}(K)-N(K')} \to Fr_{S^3-N(K)}$, this map in turn induces a pull-back map 
$$(f_*)^* :H^2(Fr_{S^3-N(K)};\mathbb{Z}) \to H^2(Fr_{Y-N(K')};\mathbb{Z})$$
Picking a trivialization of the tangent bundle of $S^3-N(K)$ and then a compatible trivialization of the tangent 
bundle of $Y-N(K')$ (the compatibility facilitated by $f$ in the obvious way), it is easy to see that the 
restriction of $(f_*)^*$ to $H^2(SO(3);\mathbb{Z})$ is injective. Therefore, $(f_*)^*$ descends to a map 
\begin{equation} \label{spincpullback}
f^*:Spin^c(S^3-N(K)) \to Spin^c(Y-N(K'))
\end{equation}
A similar discussion for spin-structure also yields a map, still denoted by $f^*$, which with the previously defined one fits into the commutative diagram
%
%
%
\begin{equation} \label{commutativediagramone}
\xymatrix{
Spin( Y-N(K'))  \ar[r]   &  Spin^c( Y-N(K')) \\
Spin( S^3-N(K))  \ar[r] \ar[u]^{f^*}  &  Spin^c( S^3-N(K)) \ar[u]_{f^*} 
}
\end{equation}
We shall return to $f^*$ after proving the next auxiliary lemma.  
\begin{lemma} 
Let $f:Y_{p^n}(K)\to S^3$ be the $p^n$-fold cyclic covering map with branching set the knot $K\subset S^3$ and set $K'=f^{-1}(K)$. Then the restriction map $H^2(Y_{p^n}(K);\mathbb{Z})\to H^2(Y_{p^n}(K)-N(K');\mathbb{Z})$ is an isomorphism. Consequently, the restriction map 
$$Spin^c(Y_{p^n}(K))\to Spin^c(Y_{p^n}(K)-N(K'))$$
is likewise an isomorphism (of affine spaces). 
\end{lemma}
\begin{proof}
Set $Y=Y_{p^n}(K)$ and consider the Mayer-Vietoris sequence in cohomology for the decomposition 
$Y = (Y-N(K')) \cup N(K')$:
\begin{align} \nonumber
0\to H^1(Y-N(K');\mathbb{Z})\oplus H^1(N(K');\mathbb{Z}) \to & H^1(\partial N(K');\mathbb{Z}) \to \cr 
\to & H^2(Y;\mathbb{Z}) \to H^2(Y-N(K');\mathbb{Z}) \to 0
\end{align}
Since $ H^1(Y-N(K');\mathbb{Z})$ is generated by the Hom-dual of the meridian of $K'$, we see that the map 
$ H^1(Y-N(K');\mathbb{Z})\oplus H^1(N(K');\mathbb{Z}) \to  H^1(\partial N(K');\mathbb{Z})$ is an isomorphism 
and therefore so is the map $H^2(Y;\mathbb{Z}) \to H^2(Y-N(K');\mathbb{Z})$.
\end{proof}
\begin{corollary} 
Every spin$^c$-structure on $Y_{p^n}(K)-N(K')$ extends, in a unique manner, to a 
spin$^c$-structure on $Y_{p^n}(K)$. 
\end{corollary}
\begin{definition} \label{ofs0}
Let $K$ be a knot in $S^3$, $p$ a prime integer and $n\ge 1$ a natural number. Let $f:Y_{p^n}(K) \to S^3$ 
be the $p^n$-fold branched covering map with branching set $K$. We define $\s_0 = \s_0(K,p^n)  \in Spin^c(Y_{p^n}(K))$ to be the unique spin$^c$-structure whose restriction to ${Y_{p^n}(K) - N(f^{-1}(K))}$ is the 
pull-back spin$^c$-structure $f^*(\s)$ (see \eqref{spincpullback}) of the unique spin$^c$-structure 
$\s\in Spin^c(S^3-N(K))$. We shall refer to $\s_0$ as the canonical spin-structure of $(K,p^n)$.
\end{definition}
\begin{theorem}  \label{s0properties}
The canonical spin-structure of a knot satisfies the following properties. 
\begin{enumerate}
\item If $K\subset S^3$ is a smoothly slice knot with slice disk $D^2 \hookrightarrow D^4$ and if $X$ is the 
the $p^n$-fold branched cover of $D^4$ with branching set $D^2$ (so that $\partial X = Y_{p^n}(K)$), then 
$\s_0(K,p^n)$ lies in the image of the restriction map $ Spin^c(X) \to Spin^c(Y_{p^n}(K))$.
\item If $K_1, K_2\subset S^3$ are two knots, then (see \eqref{addingspinc})
$$\s_0(K_1\# K_2,p^n) = \s_0(K_1,p^n) \# \s_0(K_2,p^n)$$
\item $\s_0(K,p^n)$ is a spin-structre on $Y$. 
\end{enumerate}
\end{theorem}
\begin{proof} 1. Let us denote the slice disk $D^2\hookrightarrow D^4$ by $\sigma$ and let 
$F:X\to D^4$ be the branched covering map and let $f=\partial F$. Set $\sigma ' = F^{-1}(\sigma)$ and let $N(\sigma)$ 
and $N(\sigma') = F^{-1}(N(\sigma))$ denote tubular neighborhoods of $\sigma$ and $\sigma '$ 
respectively. Note that $N(\sigma) \cong D^2 \times D^2 \cong N(\sigma ')$. 

Since $F:X-N(\sigma ') \to D^4-N(\sigma)$ is a local diffeomorphism, it induces a pull-back map 
$F^*:Spin^c(D^4-N(\sigma)) \to Spin^c(X-N(\sigma'))$. The proof of this analogous to that of $f^*$ discussed in equation \eqref{spincpullback}, the details are omitted. Of course, $Spin^c(D^4)$ is a one-point set and we denote
its sole member by $\t$. By abuse of notation, we let $\t$ also 
denote the only element of $Spin^c(D^4-N(\sigma))$. Clearly $\t|_{S^3-N(K)}=\s$, the unique 
spin$^c$-structure on $S^3-N(K)$.  Consider now the commutative diagram

\centerline{
\xymatrix{
 Spin^c(X ) \ar[d]^{r_2}  \ar[r]^{r_1}  &Spin^c(Y)  \ar[d]^{r_3}_{\cong}   \\
Spin^c(X-N(\sigma'))   \ar[r]^{r_4} &  Spin^c(Y -N(K)') \\
Spin^c(D^4-N(\sigma)) \ar[u]^{F^*}   \ar[r]^{r_5} &  Spin^c(S^3 -N(K))  \ar[u]_{f^*}\\
}
}
\vskip1mm
\noindent 
where all $r_j$ stand for restriction maps. Then 
$$\s_0(K,p^n)  = r_3^{-1}(f^*(r_5(\t))) = r_3^{-1}(r_4(F^*(\t)))$$
Thus, to prove point 1 of the theorem, we need to show that $F^*(\t) \in Im(r_2)$. 
This in turn follows from the Mayer-Vietoris sequence for the decomposition $X = (X-N(\sigma '))\cup N(\sigma ')$
by which the restriction map $H^2(X;\mathbb{Z})\to H^2(X-N(\sigma ');\mathbb{Z})$ induces an isomorphism.

\vskip2mm
2. Consider two copies $S_1$ and $S_2$ of $S^3$ with $K_i \subset S_i$. Fix identifications of $N(K_i)$ 
with $S^1\times D^2$ and pick small unknotted arcs $I_i\subset K_i$. We shall use the 3-balls 
$B_i = I_i\times D^2 \subset N(K_i)$ to perform the connected sum of $(S_1,K_1)$ and $(S_2,K_2)$, i.e. 
$$(S^3,K_1\#K_2) = ((S_1-B_1) \cup_{\partial B_1=\partial B_2} (S_2-B_2),   
(K_1-I_1)\cup_{\partial I_1=\partial I_2} (K_2-I_2))$$ 
Let $Y_i = Y_{p^n}(K_i)$ and $Y=Y_{p^n}(K_1\#K_2)$ and let 
$f_i:Y_i\to S^3$ and $g:Y\to S^3$ be the $p^n$-fold branched covering maps. Note that $g=f_1\# f_2$. 
As before, let $K_i'=f_i^{-1}(K_i)$ and define $I_i' \subset K_i'$ as $f_i^{-1}(I_i)$ and 
set $B'_i = I_i'\times D^2 = f_i^{-1}(B_i) \subset N(K_i')$. Observe that 
$Y_{p^n}(K_1\#K_2) \cong Y_{p^n}(K_1)\# Y_{p^n}(K_2)$ where the connected sum is performed by removing 
$B'_i$ from $Y_{p^n}(K_i)$ and gluing them along their boundaries, exercising care so as to glue $K_1' - f_1^{-1}(I_1)$ to $K_2'-f_2^{-1}(I_2)$ in an 
orientation respecting manner.  

The complement of the tubular neighborhood $N(K_1\#K_2)$ in $S^3$ can be obtained from  
$S^3-N(K_1)$ and $S^3-N(K_2)$ by gluing them along $I_i\times \partial B_i$, a similar statement holds for $Y_{p^n}(K_1\#K_2) - N(K_1'\#K_2')$:
\begin{align} \nonumber
S^3-N(K_1\#K_2) = (S^3-N(K_1)) \cup_{I_1\times \partial B_1 = I_2 \times \partial B_2} (S^3-N(K_2)) \cr
Y-N(K'_1\#K'_2) = (Y_1-N(K'_1)) \cup_{I'_1\times \partial B'_1 = I'_2 \times \partial B'_2} (Y_2-N(K'_2)) 
\end{align}
A Mayer-Vietoris argument now shows that the restriction maps induce isomorphisms 
\begin{align} \nonumber
H^2(S^3-N(K_1\#K_2);\mathbb{Z}) & \cong H^2(S^3-N(K_1);\mathbb{Z}) \oplus H^2(S^3-N(K_2);\mathbb{Z}) \cr
H^2(Y-N(K'_1\#K'_2);\mathbb{Z}) & \cong H^2(Y_1-N(K'_1);\mathbb{Z}) \oplus H^2(Y_2-N(K'_2);\mathbb{Z}) 
\end{align}
descending to affine isomorphisms between the corresponding spaces of spin$^c$-structures. 
This leads to the commutative diagram

\centerline{
\xymatrix{
Spin^c(Y-N(K'_1\# K'_2)) \ar[r]   &  Spin^c(Y_1-N(K'_1)) \times Spin^c(Y_2-N(K'_2))   \\
Spin^c(S^3-N(K_1\# K_2)) \ar[u]^{g^*}  \ar[r] &  Spin^c(S^3 -N(K_1)) \times Spin^c(S^3- N(K_2))  \ar[u]_{f_1^*\times f_2^*} \\
}
}
\vskip1mm
\noindent 
from which the proof of point 2 of the theorem follows. Namely, let $\s \in Spin ^c(S^3-N(K_1\#K_2))$ be the unique
spin$^c$-structure on $S^3-N(K_1\#K_2)$ so that $\s_0(K_1\#K_2,p^n)$ is the unique extension to $Y$ of 
$g^*(\s)$. But $\s_0(K_i,p^n)$ in turn is the unique extension to $Y_i$ of $f_i^*(\s)$ and so the commutativity 
of the above diagram along with $g=f_1\#f_2$, shows that $\s_0(K_1\#K_2) = \s_0(K_1,p^n)\# \s_0(K_2,p^n)$, as claimed.  

3. This is trivial, it follows for example from the commutative diagram \eqref{commutativediagramone}. 
\end{proof}
\begin{remark}
Our definition of $\s_0(K,p^n)$ agrees, as a spin$^c$-structure, with that given by Grigsy, Ruberman and Strle in lemma 2.1 of \cite{elidannysaso}. The definition from \cite{elidannysaso} takes extra care to 
extend $f^*(\s)$ to $Y_{p^n}(K)$ as a spin-structure and thus requires, in some cases, twisting by an 
element of $H^1(Y;\mathbb{Z}_2)$. In those cases their definition differs from ours as a spin-structure, though
not as a spin$^c$-structure.  
\end{remark}
%
%
%
\section{Background material} \label{backgroundmaterial}
This section gathers some supporting material to be used in section \ref{linear-independence} in the proof of theorem \ref{main}.
\subsection{Surgeries on torus knots}
We start by recalling a beautiful theorem due to Louise Moser explaining which 3-manifolds 
are obtained by surgeries on torus knots. 
\begin{theorem}[Moser \cite{moser}] \label{torusknotsurgeries}
Let $a,b \in \mathbb{N}$ be two nonzero and relatively prime integers. Let $T_{a,b} \subset S^3$ denote the $(a,b)$ torus knot and let $S^3_{p/q}(T_{a,b})$ denote the 3-manifold obtained by $p/q$-framed Dehn surgery on $T_{a,b}$.  Then 
$$S_{p/q}^3(T_{(a,b)}) = \left\{
\begin{array}{cl}
L(|p|,qb^2) & \quad \mbox{ if } \quad |abq-p|=1 \cr
L(a,b)\# L(b,a) & \quad \mbox{ if } \quad |abq-p| =0 \cr 
S(a,b,|abq-p|) &  \quad \mbox{ if } \quad |abq-p| >0 \cr
\end{array}
\right.$$
where $S(x_1,x_2,x_3)$ with $x_i\in\mathbb{Z}-\{0\}$ is a Seifert fibered space with 3 singular fibers of multiplicities $x_1,x_2$ and $x_3$. 
\end{theorem}
In section 3 of \cite{moser}, Moser outlines an algorithm to pin down the exact Seifert fibered space $S(x_1,x_2,x_3)$ in case 3 of the 
above theorem. Rather than addressing how to do this in general, we focus on a special case of 
interest to us. Recall first that the {\em Brieskorn sphere $\Sigma (a,b,c)$} associated to a triple $a,b,c\in \mathbb{N}$ of mutually prime integers, is the integral homology $3$-sphere obtained as the intersection 
$V(a,b,c) \cap S^5\subset \mathbb{C}^3$ where $V(a,b,c)$ is the complex variety
$$V(a,b,c) = \{ (z_1,z_2,z_3)\in \mathbb{C}^3 \, | \, z_1^a+z_2^b+z_3^c=0 \}$$
which is smooth away from the origin. Being the boundary of $V(a,b,c)\cap D^6$, $\Sigma (a,b,c)$ carries a natural orientation. The Brieskorn sphere $\Sigma (a,b,c)$ is realized as (see Kauffman \cite{kauffman1}):
\begin{equation} \label{branchedcoversoftorusknots}
 \Sigma (a,b,c) = \left\{ 
\begin{array}{l}
\mbox{The $a$-fold branched cover of the torus knot $T_{(b,c)}$.} \cr
\mbox{The $b$-fold branched cover of the torus knot $T_{(a,c)}$.} \cr
\mbox{The $c$-fold branched cover of the torus knot $T_{(a,b)}$.} 
\end{array}
\right.
\end{equation}
With the above orientation convention, we get 
\begin{corollary}[Moser \cite{moser}]  \label{Brieskornsurgeries}
Let $n$ be a natural number and $(a,b)$ a pair of coprime positive integers. Then $\pm1/n$ surgery on the torus knot $T_{(a,b)}$ yields the Brieskorn sphere $-\Sigma (a,b,abn\mp 1)$. 
\end{corollary}
\subsection{Heegaard Floer correction terms}
The material presented in this section can be found in \cite{peter4}. 

Let $(Y,\s)$ be a pair consisting of a rational homology 3-sphere $Y$ and a spin$^c$-structure $\s \in Spin^c(Y)$. 
To such a pair, Ozsv\'ath and Szab\'o \cite{peter4} associate the rational number $d(Y,\s) \in \mathbb{Q}$ called 
the {\em correction term} of $(Y,\s)$. The correction terms satisfy a number of properties, some of which we point to in the next theorem.
\begin{theorem}[Ozsv\'ath - Szab\'o, \cite{peter4}]  \label{dprops}
The correction terms satisfy the three properties: 
\begin{enumerate}
\item If $X$ is a rational homology $4$-ball with boundary $Y$ and if $\s$ is a spin$^c$-structure on $Y$ lying in the image of the map $Spin^c(X)\to Spin^c(Y)$, then $d(Y,\s)=0$.
\item If $(Y_1,\s_1)$ and $(Y_2,\s_2)$ are two spin$^c$ rational homology 3-spheres, then 
$$ d(Y_1\#Y_2, \s_1\#\s_2) = d(Y_1,\s_1) + d(Y_2,\s_2) $$
\item If $-Y$ denotes $Y$ with reversed orientation, then $d(-Y,\s) = -d(Y,\s)$.
\end{enumerate} 
\end{theorem}
The correction terms are in general hard to compute. In select cases however, Ozsv\'ath and Szab\'o provide 
easy to use formulae computing them. One such formula is
\begin{theorem}[Ozsv\'ath - Szab\'o, \cite{peter4}] \label{correctiontermsoftorusknots}
Let $K\subset S^3$ be a knot and suppose that $p$-framed surgery on $K$, with $p>0$, yields a lens space. For a rational number $r$ let $S^3_r(K)$ 
denote the result of $r$-framed Dehn surgery on $K$. Then, for the unique spin$^c$-structure $\s_0$ on $S_{\pm 1/n}(K)$, we obtain 
$$ d(S^3_{1/n}(K),\s_0) = -2t_0 \quad \quad \mbox{ and } \quad \quad d(S^3_{-1/n}(K),\s_0) = 0 $$ 
where $t_0$ is the $0$-th torsion coefficient computed from the symmetrized Alexander polynomial 
$\Delta_K(t)= a_0 +\sum _{j=1}^d a_i(t^i +t^{-1})$  of $K$ as $t_0 = \sum _{j=1}^d j a_{j}$. 
\end{theorem}
\begin{corollary} \label{theds}
The correction terms of $\Sigma (a,b,nab\pm 1)$ are 
\begin{equation} \label{dofsigma}
d(\Sigma(a,b,nab-1)) <0  \quad \quad \mbox{ and } \quad \quad d(\Sigma(a,b,abn+1)) = 0
\end{equation}
\end{corollary}
\begin{proof}
This corollary is a direct consequence of theorem \ref{torusknotsurgeries}, corollary \ref{Brieskornsurgeries} and theorem \ref{correctiontermsoftorusknots} along with computing $t_0$ for torus knots. Since the Alexander polynomial $\Delta_{T_{(a,b)}}(t)$ of the torus knot $T_{(a,b)}$ equals 
$$\Delta_{T_{(a,b)}}(t) = \frac{(t^{ab}-1)(t-1)}{(t^a-1)(t^b-1)}  \cdot t^{-(a-1)(b-1)/2}$$
its coefficients are $\pm 1$ and alternate in sign. The leading coefficient of $\Delta_{T_{(a,b)}}(t) $ is 1
rendering the $0$-th torsion coefficient $t_0$ of $T_{(a,b)}$ positive and therefore making $-2t_0$ negative, as claimed. 
\end{proof}
\section{Linear independence} \label{linear-independence}
Recall from the introduction that we defined $\delta_{p^n}(K)$ to be $2 d(Y_{p^n}(K),\s_0(K,p^n))$ with 
$\s_0(K,p^n)$ as specified in definition \ref{ofs0}. We first show that $\delta_{p^n}$ gives rise to a well defined homomorphism 
$\delta _{p^n} :\mathcal C \to \mathbb{Z}$. The proofs of the proceeding two propositions heavily rely on Kauffman's results from \cite{kauffman1, kauffman2}. 
\begin{proposition} \label{kauffmansresults}
Let $K\subset S^3$ be a knot and $\tilde \Sigma\subset S^3$ be any Seifert surface of $K$.  Let $\Sigma\subset D^4$ be obtained from $\tilde \Sigma$ by pushing the interior of the latter into $D^4$ (so that $\Sigma$ is properly embedded in $D^4$ and $\partial \Sigma = K$). Let $X_\Sigma$ be the $p^n$-fold branched cover of $D^4$ with branching set $\Sigma$. Then the following hold: 
\begin{enumerate}
\item The signature $\sigma (X_\Sigma)$ of $X_\Sigma$ is 
$$\sigma(X_\Sigma ) =  \sum_{i=0}^{p^n-1}\sigma_{\omega ^i}(K) $$
where $\sigma _{\tau}(K)$ is the Tristram-Levine signature of $K$ associated to $\tau \in S^1$ and $\omega$ is a primitive
$p^n$-th root of unity. 
\item $X_\Sigma$ is a spin manifold and with a unique spin-structure. 
\end{enumerate}
\end{proposition}
\begin{proof}
This proposition is largely contained in the work of  Kauffman \cite{kauffman1} (see also \cite{kauffman2}) 
where he extensively studies the algebraic topology of the manifold $X_\Sigma$. Specifically, the signature formula from the first claim of the proposition has been worked 
out by Kauffman (page 290 in \cite{kauffman2}). 

For the second claim we also rely on \cite{kauffman1, kauffman2}. In these works, Kauffman finds an explicit
matrix representative for the intersection form on $H^2(X_\Sigma;\mathbb{Z})$ (in terms of the linking matrix of $K$ associated to $\Sigma$, see page 283 of \cite{kauffman2}) from which one can readily pin down the second Stiefel-Whitney class 
$w_2(TX_\Sigma)$ and finds the latter to be zero implying that $X_\Sigma$ is a spin manifold. Since 
$Spin(X)$ is an affine space on $H^1(X_\Sigma;\mathbb{Z}_2)$, the second claim of the proposition follows
from $H^1(X_\Sigma;\mathbb{Z}_2)=0$, another results from \cite{kauffman2} (page 282). 
\end{proof}
\begin{proposition}
Let $p$ be any prime number. Then 
\begin{enumerate}
\item For every knot $K\subset S^3$ the number $\delta _{p^n}(K)$ is an integer. 
\item For any two knots $K_0,K_1\subset S^3$ one obtains $\delta _{p^n}(K_0\#K_1) = \delta _{p^n}(K_0) + \delta _{p^n}(K_1)$.
\item If $K$ is smoothly slice then $\delta_{p^n}(K)=0$ for all choices of $p,n\in \mathbb{N}$ with $p$ prime.  
\end{enumerate}
In particular, $K\mapsto \delta _{p^n}(K)$ is a group homomorphism from $\mathcal C \to \mathbb{Z}$. 
\end{proposition}
\begin{proof}
For the first statement of the proposition we recall a formula proved by Ozsv\'ath and Szab\'o in \cite{peter4}. To 
state the result, let $X$ be any smooth 4-manifold with $\partial X=Y_{p^n}(K)$ and let $\t \in Spin^c(X)$ be any spin$^c$-structure with  $\t |_{Y_{p^n}(K)}=\s_0(K,p^n)$. Then  
\begin{align} \nonumber
d(Y_{p^n}(K),\s_0(K,p^n)) & \equiv \frac{c_1(\t)^2 - \sigma }{4} \, \,  (\mbox{mod } 2) 
\end{align}
Given a Seifert surface $\Sigma \subset S^3$ of $K$, let $X=X_\Sigma$ be the 4-manifold from proposition 
\ref{kauffmansresults} and note that $\partial X = Y_{p^n}(K)$.  Let $F:X_\Sigma \to D^4$ be the branched covering map and let $f:Y_{p^n}(K) \to S^3$ be $\partial F$.  The restriction map $Spin^c(X)\to Spin^c(Y_{p^n}(K))$ is 
modeled on the map $H^2(X;\mathbb{Z})\to H^2(Y_{p^n}(K);\mathbb{Z})$. The latter is surjective since 
$H^3(X,Y_{p^n}(K);\mathbb{Z})\cong H_1(X;\mathbb{Z})=0$ by a result of Kauffman's (page 282 in \cite{kauffman2}). Thus every spin$^c$-structure on $Y_{p^n}(K)$ extends to a spin$^c$-structure on $X$. 

To see that $\s_0(K,p^n)$ extends to the unique spin-structure $\t_0$ on $X$ (part 2 of proposition 
\ref{kauffmansresults}), consider the following isomorphisms
$$ H^2(X,X-N(\Sigma') ;\mathbb{Z}_2) \cong H^2(N(\Sigma '),\partial N(\Sigma ');\mathbb{Z}_2) \cong 
H_2(N(\Sigma ');\mathbb{Z}_2) \cong 0 $$
The first of these follows by excision and the second by Alexander-Poincar\'e duality. With this as input, consider
the following portion of the exact sequence of the pair $(X,X-N(\Sigma '))$ with $\mathbb{Z}_2$-coefficients:
$$ ... \to H^1(X;\mathbb{Z}_2) \to H^1(X-N(\Sigma');\mathbb{Z}_2) \to H^2(X,X-N(\Sigma');\mathbb{Z}_2) \to ...$$
As already pointed out in the proof of proposition \ref{kauffmansresults}, Kauffman's results from \cite{kauffman2} show that $H^1(X;\mathbb{Z}_2)=0$ and so we conclude that $H^1(X-N(\Sigma ');\mathbb{Z}_2)=0$ also. This 
shows that $Spin(X-N(\Sigma '))$ consists of a single spin-structure and that therefore the restriction map 
$Spin(X) \to Spin(X-N(\Sigma '))$ is an isomorphism. The fact that $\s_0(K,p^n)$ lies in the image of 
$Spin(X) \to Spin(Y_{p^n}(K))$ now follows from the commutative diagram (with all $r_j$ being restriction maps):

\centerline{
\xymatrix{
 Spin^c(X ) \ar[d]^{r_2}  \ar[r]^{r_1}  &Spin^c(Y)  \ar[d]^{r_3}_{\cong}   \\
Spin^c(X-N(\Sigma'))   \ar[r]^{r_4} &  Spin^c(Y -N(K)') \\
Spin^c(D^4-N(\Sigma)) \ar[u]^{F^*}   \ar[r]^{r_5} &  Spin^c(S^3 -N(K))  \ar[u]_{f^*}\\
}
}
\noindent Namely, $\s_0(K,p^n) = r_3^{-1}(f^*(r_5(\t))=r_3^{-1}(F^*(r_4(\t))$ where $\t\in Spin^c(D^4-N(\Sigma))$ 
is the unique spin-structure on $D^4-N(\Sigma)$. Since $F^*(\t)$ is also a spin-structure it must extend to the (unique) spin-structure $\t_0$ on $X$ and so $s_0(K,p^n) = r_1(\t_0)$ as claimed. 

Finally, since $\delta_{p^n}(K) = 2 d(Y_{p^n}(K).\s_0(K,p^n))$ and the latter is congruent to $(c_1(\t_0)^2-\sigma)/2 = -\sigma /2$ modulo $2$, we see that it must be an integer. 
\vskip1mm
The second and third statement of the proposition follow directly from theorems \ref{s0properties} and \ref{dprops}. 
\vskip1mm
With statements 1--3 proved, it is now automatic for $\delta_{p^n}$ to descend to a group 
homomorphism $\delta_{p^n}:\mathcal C \to \mathbb{Z}$. First of, note that $\delta_{p^n}(-K) = -\delta_{p^n}(K)$ (where $-K$ is the reverse mirror of $K$) as follows from points $2$ and $3$ of the proposition along with the fact that $K\#(-K)$ is smoothly slice. If $K_1$ is smoothly concordant to $K_2$ then $K_1\#(-K_2)$ is smoothly slice so that
$$0 = \delta_{p^n}(K_1\#(-K_2)) = \delta_{p^n}(K_1) + \delta_{p^n}(-K_2) = \delta_{p^n}(K_1) - \delta_{p^n}(K_2)$$
\end{proof}
\begin{definition} \label{fermatmersenne}
A natural number $n$ is called a Fermat number if it is of the form $n=2^{2^k}+1$ for some $k\in \mathbb{N}\cup \{0\}$. Similarly, $n$ is called a Mersenne number if it looks like $n=2^m-1$ with $m\in \mathbb{N}$. A natural number $n$ is called a Fermat prime or a Fermat prime power (Mersenne prime or Mersenne prime power) if it is a Fermat number (Mersenne number) and prime or a prime power at the same time. 
\end{definition}
It is well known that if a Mesenne number $2^m-1$ is a prime power, then $m$ must be a prime number itself \cite{lighneal}. It is unknown which Fermat numbers are either prime or prime powers. The only Fermat primes known to date are the fist five, namely $3$, $5$, $17$, $257$ and $65537$. 
\begin{lemma} \label{fermatmersenneprimepowers}
Let $p$ be an odd prime. Then 
\begin{enumerate}
\item $p^{2^n}-1$ can only be a prime power if $p$ is a Fermat prime. If $p=3$ this happens for $n=0,1$ and if $p>3$ this can happen only for $n=0$.
\item $p^{2^n}+1$ can only be a prime power if $p^{2^n}$ is a Mersenne prime power. With $p$ fixed, this can happen for at most one $n\in \mathbb{N}\cup \{0\}$.  
\end{enumerate} 
\end{lemma}
\begin{proof}
Pick and fix a prime $p\ge  3$ throughout the proof. Consider first when $p^{2^n}-1$ is of the form $2^m$ for some integer $m$. 
When $n$ is at least $1$, $p^{2^n}-1$ factors as  
$$p^{2^n}-1 = (p^{2^{n-1}}-1)(p^{2^{n-1}}+1)$$
and so, in order for $p^{2^n}-1$ to be a prime power, both $p^{2^{n-1}}\pm 1$ 
have to be powers of two. An easy analysis shows that this only happens if $p=3$ and $n=1$. 
If $n=0$ and $p^{2^n}-1=p-1$ is a prime power, say $2^m$, then $p=2^m+1$. Elementary number theory then 
shows that $m$ itself has to be a power of $2$ (this observation was already made by Gauss) so that $p$  becomes a Fermat prime. Since $p=3$ is itself a Fermat prime, the first statement of the lemma follows. 
\vskip1mm
Turning to $p^{2^n}+1$, suppose we found some $n$ for which $p^{2^n}+1=2^m$. Then 
$$p^{2^{n+1}}+1 = \left( p^{2^n} \right)^2 +1 = (2^m-1)^2+1  = 2(2^{2m-1}-2^m +1) $$
Since the second factor in the last term above is odd, we see that $p^{2^{n+1}}+1$ cannot be a prime power. 
A similar argument shows the same to be true for $p^{2^{n+k}}+1$ for any $k\ge 1$:
\begin{align} \nonumber
p^{2^{n+k}}+1 & = \left( p^{2^n} \right)^{2^k} +1 = (2^m-1)^{2^k}+1  = 1+ \sum _{j=0}^{2^k} {2^k \choose j} (-1)^j 2^{mj} \cr
& = 2\left( 1 + \sum _{j=1}^{2^k} {2^k \choose j} (-1)^j 2^{mj-1}\right)  = 2(1+\mbox{even number})
\end{align}
This shows that $p^{2^n}+1$ can be a prime power for at most one $n\in \mathbb{N}\cup \{ 0 \}$. 
Note that $p^{2^n}+1=2^m$ can only happen if $p^{2^n}=2^m-1$, i.e. if $p^{2^n}$ is a Mersenne prime power. 
For this to happen $m$ must itself be a prime number \cite{lighneal} . The choices of $p=7$, $n=0$ and $m=3$ are an example.
\end{proof}
The next theorem is a slightly strengthened version of theorem \ref{main}.
\begin{theorem} \label{mainprime}
Let $p\ge 3$ be a prime integer. If $p=3$ set $\mathcal F=\{0,1\}$, if $p\ne 3$ is a Fermat prime let $\mathcal F=\{0\}$ and otherwise let $\mathcal F=\emptyset$. If there exists an integer $n_0\in \mathbb{N}\cup\{0\}$ such that 
$p^{2^{n_0}}$ is a Mersenne prime power (see definition \ref{fermatmersenne}), let $\mathcal M=\{n_0\}$, otherwise let $\mathcal M=\emptyset$. Then the set 
$$ \{ \delta _{p^{2^n}}:\mathcal C \to \mathbb{Z}\, | \, n\in (\mathbb{N}\cup\{0\})  - (\mathcal F \cup \mathcal M) \} \subset Hom(\mathcal C, \mathbb{Z})$$
is linearly independent.
\end{theorem}
\begin{proof}
According to lemma \ref{fermatmersenneprimepowers}, no element of the infinite set 
$$ \{ p^{2^n}\pm 1 \, | \, n\in   (\mathbb{N}\cup\{0\})  - (\mathcal F \cup \mathcal M) \}$$
is a prime power. 
Let $n_i$, $i\in \mathbb{N}$  be the sequence enumerating (in increasing order) the elements of  
$ (\mathbb{N}\cup\{0\})  - (\mathcal F \cup \mathcal M)$. For simplicity of notation we shall write $m_i$ for $p^{2^{n_i}}$. Thus we need to prove that the set $ \{ \delta_{m_i} : \mathcal C \to \mathbb{Z} \, | \, i \in \mathbb N \} $ is is linearly independent. To see this, suppose that some linear combination of $\delta _{m_i}$'s results in the zero 
homomorphism:
\begin{equation} \label{lineardependence}
\lambda _1 \delta _{m_1} + \lambda _2 \delta _{m_2} + ...\lambda _\ell \delta _{m_\ell}  = 0
\end{equation}
Since $m_1 + 1$ is not a prime power (lemma \ref{fermatmersenneprimepowers}), we can find a pair of positive coprime integers $(a,b)$ such that 
$$  m_1 = kab-1$$
According to corollary \ref{Brieskornsurgeries}, $1/k$-framed Dehn surgery on the torus knot $T_{(a,b)}$ gives the Brieskorn sphere $-\Sigma(a,b,kab-1)$.  Since $kab-1=m_1=p^{2^{n_1}}$ we see that $\Sigma (a,b,m_1)$ is the $p^{2^{n_1}}$-fold branched cover of $T_{(a,b)}$, cf. equation \eqref{branchedcoversoftorusknots}.  
But according to corollary \ref{theds} (with the help of statement 3 from theorem \ref{dprops} to address the change of orientation) we then know that 
\begin{equation} \label{aux1ds}
 \delta _{m_1}(T_{(a,b)}) < 0
\end{equation}  
On the other hand, note that $m_1-1$ is a factor of $m_s-1$ for all $s\ge 1$ since, if $n_s=n_1+d$, then 
\begin{align} \nonumber
m_s-1 & = (m_1)^{2^d}-1 =\left( (m_1)^{2^{d-1}}-1\right) \left( (m_1)^{2^{d-1}}+1\right) \cr
& = \left( (m_1)^{2^{d-2}}-1\right) \left( (m_1)^{2^{d-2}}+1\right)  \left( (m_1)^{2^{d-1}}+1\right) \cr
&\, \,  \, \, \vdots \cr
& = (m_1-1)\cdot  \prod _{j=0}^{d-1} \left( (m_1)^{2^{j}}+1\right)
\end{align}
Thus, for all $s\ge 2$ we get $m_s = \ell_scd+1$ for some integer $\ell_s$. With this, another use of corollary \ref{theds} show that 
\begin{equation} \label{aux2ds} 
\delta _{m_s}(T_{(a,b)}) = 0 \quad \quad \forall s\ge 2
\end{equation}
Plugging $T_{(a,b)}$ into equation \eqref{lineardependence}, with the help of equatioins \eqref{aux1ds} and \eqref{aux2ds}, results in $\lambda _1=0$. 

From here on one just repeats this argument $\ell -1$ more times to obtain $\lambda _i=0$ for all $i=1,...,\ell$, we omit further details.
\end{proof}
\begin{corollary}
Let $p$ be a prime integer and $n_i$, $i\in \mathbb{N}$ be the sequence enumerating the elements 
in $(\mathbb N \cup \{0\})-(\mathcal F \cup \mathcal M)$ from theorem \ref{mainprime}. Then for 
every index $i$ there is a pair of coprime positive integers $(a_i,b_i)$ such that 
$$ \delta _{p^{2^{n_j}}}(T_{(a_i,b_i)}) =0  \quad \forall j>i \quad \quad \mbox{ and } \quad \quad \delta _{p^{2^{n_i}}}(T_{(a_i,b_i)})<0$$ 
\end{corollary}
\section{Examples part 1 - Brieskorn spheres} \label{examples1-section}
In this section, the first of three, we start to address the question of how to evaluate the invariants 
$\delta_{p^n}(K)$ for concrete knots $K$ and concrete choices of $p$ and $n$. While this is a difficult task in general, we explore a number of techniques that are successful in such computations for a substantial set of examples of $K$, $p$ and $n$. The techniques stated in this section have already been exploited for the 
proof of theorem \ref{mainprime} and are only listed for emphasis. 

Let $a$ and $b$ be two positive and relatively prime integers and consider the torus knot $T_{(a,b)}$. 
Let 
$$\Delta_{T_{(a,b)}}(t) = \frac{(t^{ab}-1)(t-1)}{(t^a-1)(t^b-1)}  \cdot t^{-(a-1)(b-1)/2}$$
be its symmetrized Alexander polynomial. Let us write $\Delta_{T_{(a,b)}}(t)$ in the form 
$\Delta_{T_{(a,b)}}(t) = a_0 + f_{(a,b)}(t) + f_{(a,b)}(t^{-1})$ with $f_{(a,b)}(t) \in t\cdot \mathbb{Z}[t]$ so that the $0$-th torsion coefficient $t_0$ of $T_{(a,b)}$ takes the form $t_0(T_{(a,b)}) = f_{(a,b)}'(1)$. 

Since $1/n$ surgery on $T_{(a,b)}$ yields the Brieskorn sphere $-\Sigma (a,b,abn-1)$ and $-1/n$ surgery on 
$T_{(a,b)}$ gives $-\Sigma(a,b,abn+1)$ (cf. corollary \ref{Brieskornsurgeries}), theorems \ref{torusknotsurgeries} and \ref{correctiontermsoftorusknots} imply that 
\begin{equation} \label{formulaone}
\delta_{p^n}(T_{(a,b)}) =  \left\{ 
\begin{array}{ll}
2f'(1) & \quad ; \quad \mbox{whenever  $abm-1=p^n$ for some } m\in \mathbb{Z} \cr
& \cr
 0 & \quad ; \quad \mbox{whenever  $abm+1=p^n$ for some } m\in \mathbb{Z} 
\end{array} \right. 
\end{equation}
In the above, $p$ is always assumed to be a prime number. 
For example, taking $a=7$ and $b=9$ we obtain 
$$
f_{(7,9)}(t) =  t -t^2+t^3-t^5+t^6-t^7+t^8-t^9+t^{10} -t^{14}+t^{15}-t^{16} +t^{17}-t^{23} +t^{24}
$$
so that $f_{(7,9)}'(1) = 8$. Therefore we obtain, for instance, 
$$\begin{array}{c} 
\delta_{5^3}(T_{(7,9)}) = 16, \,\, \delta_{251}(T_{(7,9)}) = 16,  \,\, \delta_{503}(T_{(7,9)}) = 16 ,  \,\, \delta_{881}(T_{(7,9)}) = 16\cr
\mbox{ and } \cr
\delta_{2^6}(T_{(7,9)}) = 0, \,\, \delta_{127}(T_{(7,9)}) = 0,  \,\, \delta_{379}(T_{(7,9)}) = 0 ,  \,\, \delta_{631}(T_{(7,9)}) = 0,  \,\, \delta_{757}(T_{(7,9)}) = 0.
\end{array}
$$
\section{Examples part 2 - Negative definite plumbings} \label{examples2-section}
The computational methods from section \ref{examples1-section} exclude many of the branched covers
of torus knots. For instance, when $K=T_{(2,5)}$, formula \eqref{formulaone} allows for a computation of 
$\delta_{9}(K)$ and $\delta_{11}(K)$ but not for example $\delta_3(K)$ or $\delta_7(K)$. Note that the 
$3$-fold and $7$-fold branched covers of $T_{(2,5)}$ are still Brieskorn spheres, namely $\Sigma (2,3,5)$ and $\Sigma (2,5,7)$ respectively.  

Many Brieskorn spheres $\Sigma (a,b,c)$ bound negative definite plumbings, $\Sigma (2,3,5)$ and $\Sigma (2,5,7)$ are examples. When this happens, Ozsv\'ath and Szab\'o \cite{peter6} and A. Nemethi \cite{nemethi} provide formulae for computing the correction terms of $\Sigma (a,b,c)$, provided an additional restriction on the plumbing graph yielding $\Sigma (a,b,c)$ is met. Namely, given a weighted graph $G$, we shall call a vertex $x$ of $G$ 
a {\em bad vertex} if its valence $v(x)$ (the number of edges emanating from $x$) and weight $d(x)$ satisfy the inequality $d(x) < -v(x)$. The formulae from \cite{peter6} provide a completely combinatorial algorithm for the 
computation of $d(\Sigma(a,b,c))$ for all those $\Sigma (a,b,c)$ which bound negative definite plumbing graphs with no more than two bad vertices. This algorithm, besides having been explained in \cite{peter6}, has been outlined in a number of articles, see for example \cite{jabukanaik} and \cite{greenejabuka}. We thus omit it here and instead focus on an example. We would like to point out that the results from examples \ref{ex2} -- \ref{ex8} not covered by formula \eqref{formulaone}, have been computed in this manner. The results from example \ref{ex1} are addressed in the next section.

\begin{figure}[htb!] 
\centering
\includegraphics[width=14cm]{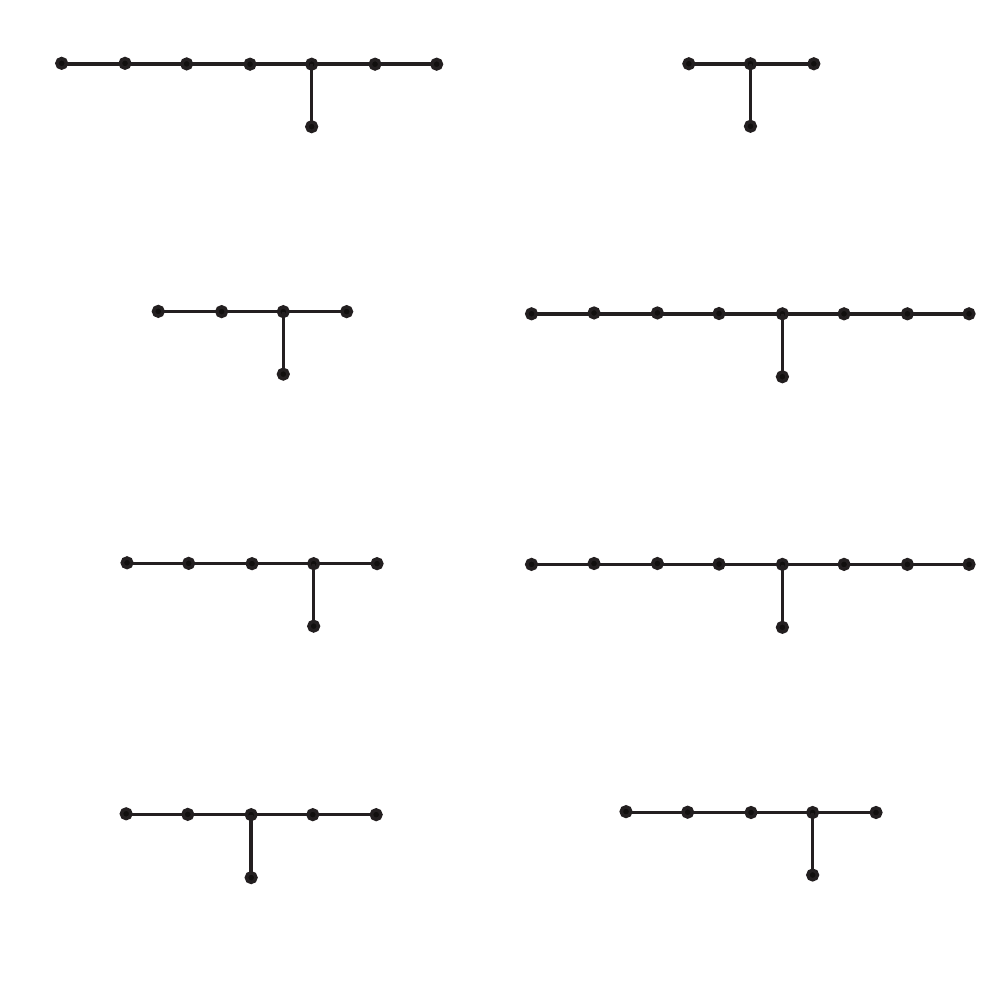}
\put(-325,320){$\Sigma(2,3,5)$}
\put(-386,380){$-2$}
\put(-361,380){$-2$}
\put(-336,380){$-2$}
\put(-311,380){$-2$}
\put(-286,380){$-2$}
\put(-261,380){$-2$}
\put(-236,380){$-2$}
\put(-295,346){$-2$}
\put(-125,320){$\Sigma(2,3,7)$}
\put(-137,380){$-2$}
\put(-112,380){$-1$}
\put(-87,380){$-3$}
\put(-120,346){$-7$}
\put(-325,220){$\Sigma(2,5,7)$}
\put(-348,280){$-2$}
\put(-323,280){$-4$}
\put(-298,280){$-1$}
\put(-273,280){$-2$}
\put(-305,246){$-5$}
\put(-125,220){$\Sigma(2,5,13)$}
\put(-199,280){$-2$}
\put(-174,280){$-2$}
\put(-149,280){$-2$}
\put(-124,280){$-2$}
\put(-99,280){$-2$}
\put(-74,280){$-2$}
\put(-49,280){$-2$}
\put(-24,280){$-5$}
\put(-106,246){$-7$}
\put(-325,120){$\Sigma(2,5,17)$}
\put(-360,180){$-3$}
\put(-335,180){$-2$}
\put(-310,180){$-4$}
\put(-285,180){$-1$}
\put(-260,180){$-2$}
\put(-293,146){$-5$}
\put(-125,120){$\Sigma(2,7,9)$}
\put(-199,180){$-3$}
\put(-174,180){$-2$}
\put(-149,180){$-2$}
\put(-124,180){$-2$}
\put(-99,180){$-2$}
\put(-74,180){$-2$}
\put(-49,180){$-2$}
\put(-24,180){$-3$}
\put(-106,146){$-2$}
\put(-325,20){$\Sigma(2,9,11)$}
\put(-360,80){$-2$}
\put(-335,80){$-5$}
\put(-310,80){$-1$}
\put(-285,80){$-4$}
\put(-260,80){$-3$}
\put(-318,46){$-2$}
\put(-125,20){$\Sigma(2,9,13)$}
\put(-162,80){$-2$}
\put(-137,80){$-3$}
\put(-112,80){$-3$}
\put(-87,80){$-1$}
\put(-62,80){$-9$}
\put(-106,46){$-2$}
\caption{Negative definite plumbing descriptions of some of the Brieskorn spheres that show up in examples \ref{ex2} -- \ref{ex8}. Note that no weighted graph has more than two bad vertices. }  \label{pic1}
\end{figure}
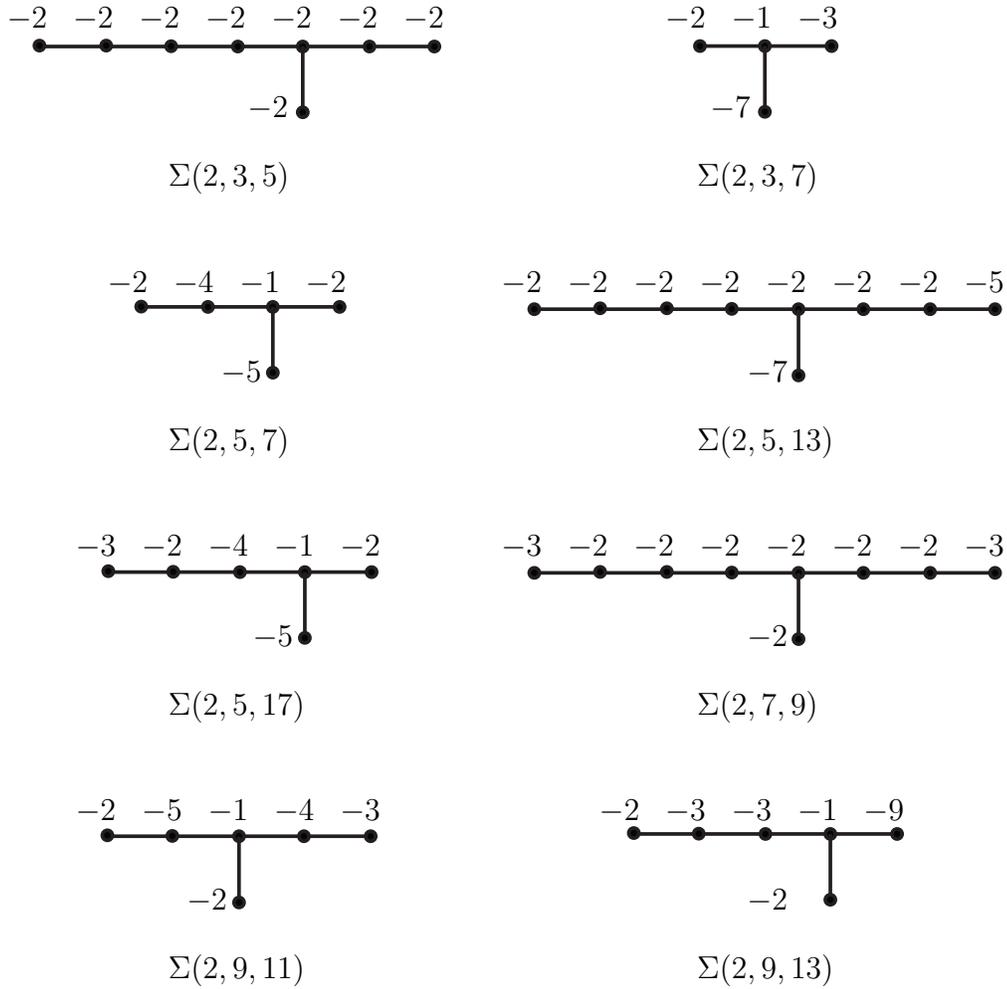

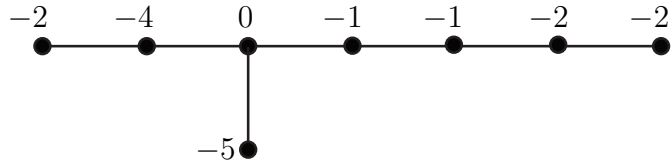
\begin{figure}[htb!] 
\centering
\includegraphics[width=10cm]{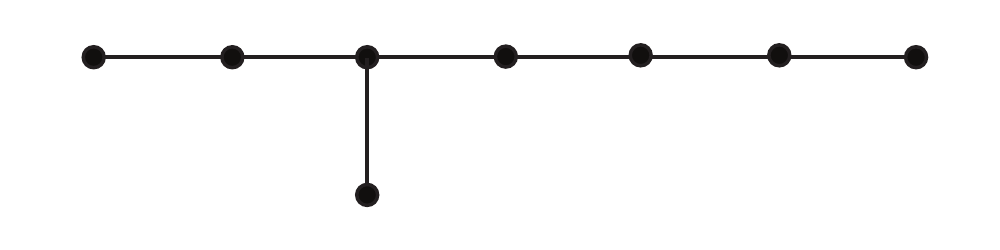}
\put(-271,62){$-2$}
\put(-231,62){$-4$}
\put(-184,62){$0$}
\put(-152,62){$-1$}
\put(-114,62){$-1$}
\put(-74,62){$-2$}
\put(-36,62){$-2$}
\put(-200,12){$-5$}
\caption{A plumbing description of the Brieskorn sphere $\Sigma(2,5,7)$. }  \label{pic2}
\end{figure}

Some of the negative definite plumbings with two or fewer bad vertices, utilized in examples \ref{ex2} -- \ref{ex8}, are described in figure \ref{pic1}. The remaining cases are left as an easy exercise. To obtain such plumbings we 
have followed the algorithm outlined in \cite{nik1,nik2} where we refer the interested reader for full details. 
By way of example, consider $\Sigma(2,5,7)$ which is the $7$-fold branched cover of $T_{(2,5)}$ (and of course the double branched cover of $T_{(5,7)}$ and the $5$-fold cover of $T_{(2,7)}$). Set $a_1=2$, $s_2=5$ and $a_3=7$ and solve the equation
$$ a_1a_2b_1+a_1b_2a_3+ b_1a_2a_3 = 1$$
for $b_1$, $b_2$ and $b_3$. For example $b_1=1$, $b_2=-1$ and $b_3=-2$ will do. Find continued 
fraction expansions for $a_i/b_i$ next:
$$\frac{a_1}{b_1} = [ -1,-1,-2,-2] \quad \quad \quad \frac{a_2}{b_2} = [-5] \quad \quad \quad 
\frac{a_3}{b_3} = [-4,-2]$$
where by $[x_1,x_2,...,x_n]$ we mean  
$$ [x_1,x_2,...,x_n] =   x_1-\cfrac{1}{x_2-\cfrac{1}{\ddots-\cfrac{1}{x_n}}}. $$
The plumbing diagram describing $\Sigma (2,5,7)$ is then obtained by drawing a single 3-valent vertex 
with framing zero from which 3 branches emerge whose vertices have framing coefficients given 
by the above 3 continued fraction expansion coefficients. This gives the plumbing diagram in figure \ref{pic2}.
Repeatedly blowing down the $-1$ framed vertices from figure \ref{pic2} leads to the plumbing description for $\Sigma (2,5,7)$ from figure \ref{pic1}. This weighted graph has only one bad vertex and its incidence matrix 
$$\left[
\begin{array}{rrrrr}
-2  &  1   &  0  &   0 &   0 \cr
1   &  -4  &  1  &   0 &   0 \cr
0   &   1  &  -1 &   1 &  1 \cr
0   &   0   &  1 &   -2 & 0 \cr
0   &   0   &   1&   0 & -5  
\end{array}
\right]$$
is easily verified to be negative definite (see for example \cite{gilbert}). 
The Ozsv\'ath-Szab\'o algorithm from \cite{peter6} now allows for a computation of 
$\delta_{7}(T_{(2,5)}) = \delta_2(T_{(5,7)}) = \delta _5 (T_{(2,7)})$ and yields $0$. We performed the needed computations by computer.  
\section{Examples part 3 - Surgery descriptions of branched covers} \label{trefoil-section}
This section computes $\delta_{2^n}(K)$ for all $n\in \mathbb{N}$ where $K$ is the right-handed trefoil. The 
results thus obtained are those from example \ref{ex1}. 
We compute these invariants by first finding an explicit Dehn surgery description of $Y_{2^n}(K)$. By happenstance, each 
$Y_{2^n}(K)$ turns out to be surgery on a single alternating knot, the knot itself depending on $n$. 
Such a description of $Y_{2^n}(K)$ allows us to use the existing Heegaard Floer 
machinery \cite{peter7} to compute the correction terms $d(Y_{2^n}(K),\s_0(K,2^n))$. 

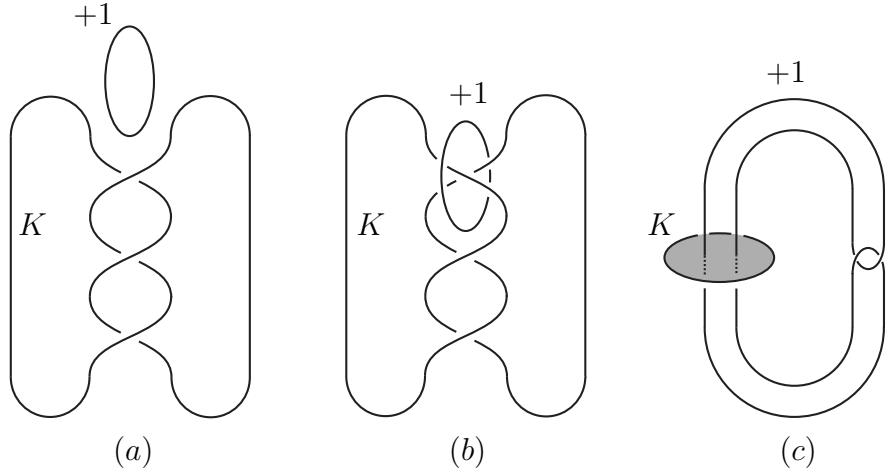
\begin{figure}[htb!] 
\centering
\includegraphics[width=12cm]{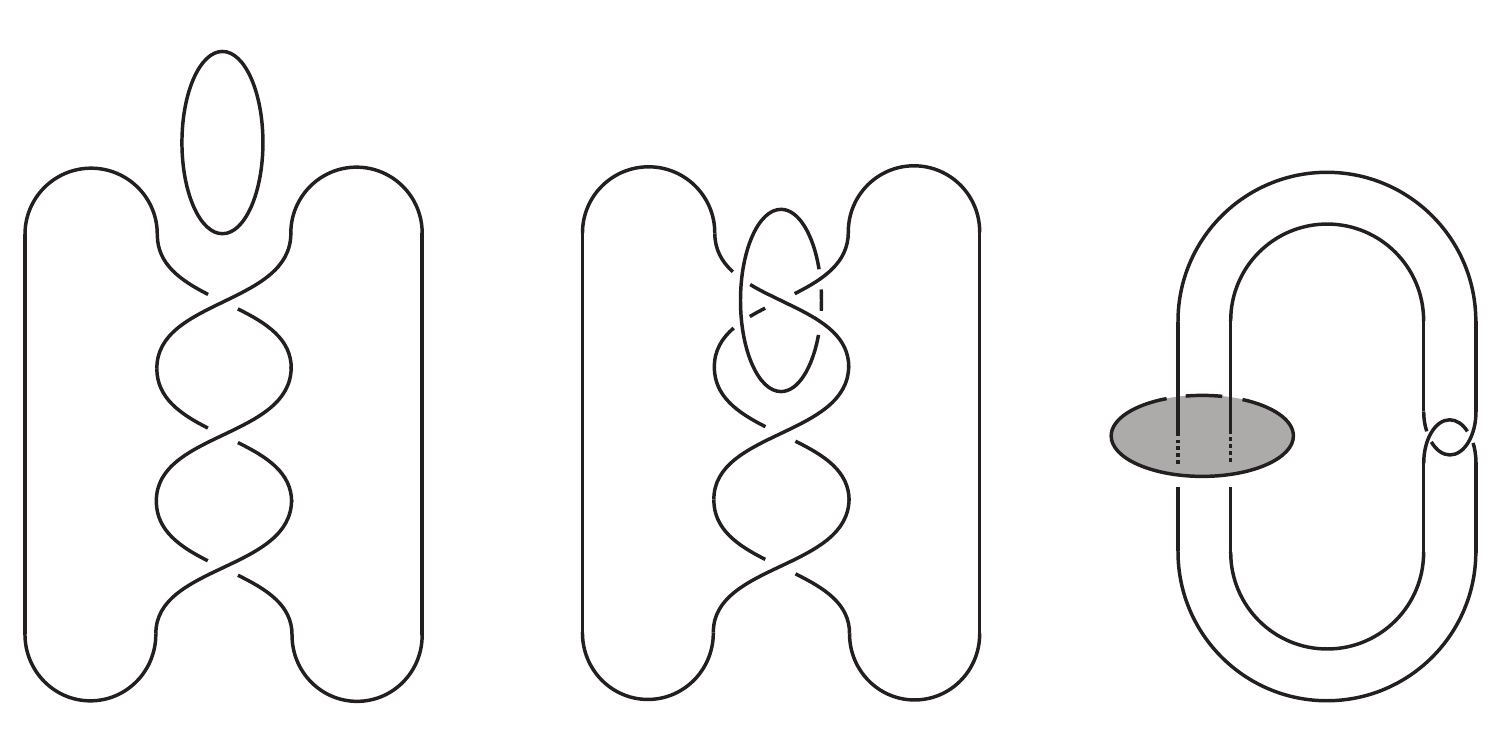}
\put(-312,160){$+1$}
\put(-170,130){$+1$}
\put(-50,138){$+1$}
\put(-333,80){$K$}
\put(-205,80){$K$}
\put(-95,80){$K$}
\put(-297,-5){$(a)$}
\put(-170,-5){$(b)$}
\put(-45,-5){$(c)$}
\caption{The right-handed trefoil $K$ represented as the unknot in $S^3$ when the latter is 
viewed as $+1$ surgery on another unknot. The disk $D$ used in the construction of $Y_m(K)$ is 
shaded in the right-most figure.}  \label{pic3}
\end{figure}

To arrive at a surgery description of $Y_{2^n}(K)$, note that the right-handed trefoil $K$ can be described as 
the unknot in $S^3$ if $S^3$ is viewed as $+1$ surgery on another unknot as in figure \ref{pic3}.
The $m$-fold cover $Y_m(K)$ of $S^3$ branched along $K$ can then be constructed by cutting $S^3$ open along the  
open disk $D$ bounded by $K$ (figure \ref{pic3}c), taking $m$ copies of $S^3-D$ and gluing the $i$-th copy to the 
$(i+1)$-st copy along their boundary components as prescribed by the orientation (and likewise gluing the $m$-th copy of $S^3-D$ to the $1$-st copy), see \cite{rolfsen} for full details. 
\begin{figure}[htb!] 
\centering
\includegraphics[width=14cm]{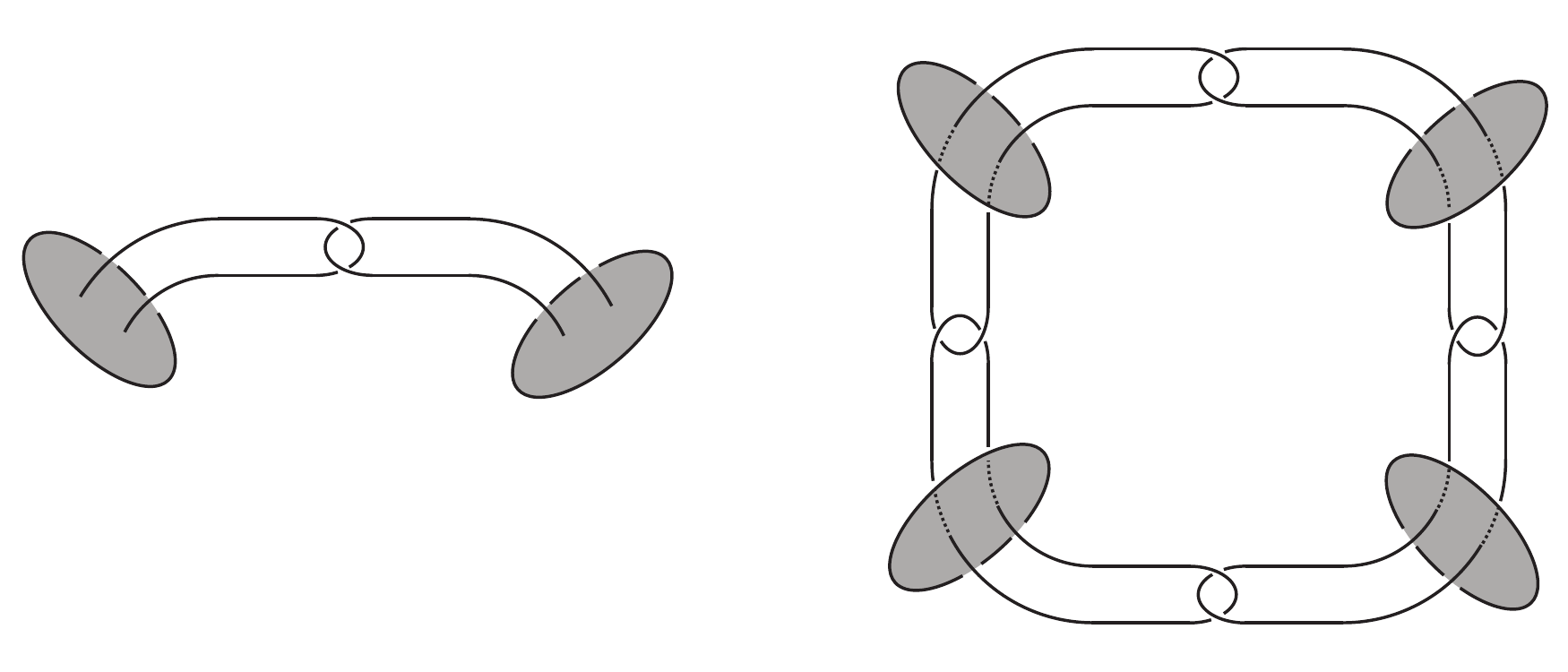}
\put(-315,-5){$(a)$}
\put(-95,-5){$(b)$}
\put(-125,130){$-1$}
\put(-75,130){$-1$}
\put(-125,33){$-1$}
\put(-75,33){$-1$}
\caption{The left figure is obtained from figure \ref{pic3}c by cutting open $S^3$ along the disk $D$. The right figure shows $4$ concatenated copies of the left figure, thus giving a surgery description of $Y_4(K)$.  }  \label{pic6}
\end{figure}
\begin{figure}[p] 
\centering
\includegraphics[width=16cm]{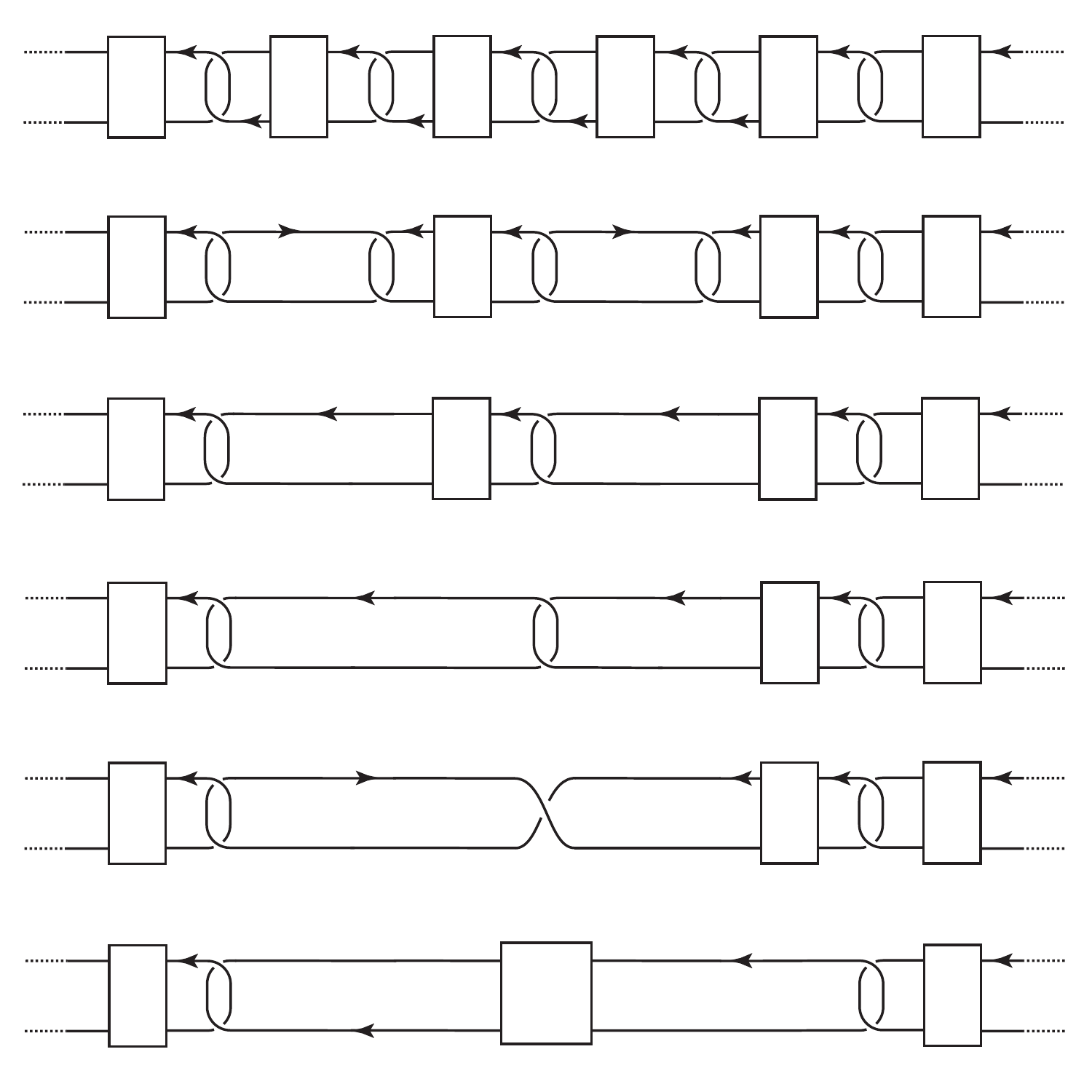}
\put(-402,415){$m$}
\put(-334,415){$n$}
\put(-265,415){$n$}
\put(-198,415){$n$}
\put(-129,415){$n$}
\put(-60,415){$k$}
\put(-402,340){$m$}
\put(-268,340){$2n$}
\put(-132,340){$2n$}
\put(-60,340){$k$}
\put(-402,265){$m$}
\put(-268,265){$2n$}
\put(-132,265){$2n$}
\put(-60,265){$k$}
\put(-345,287){$1$}
\put(-200,287){$0$}
\put(-423,287){$0$}
\put(-402,188){$m$}
\put(-132,188){$4n$}
\put(-60,188){$k$}
\put(-345,210){$1$}
\put(-200,210){$0$}
\put(-423,210){$0$}
\put(-402,114){$m$}
\put(-132,114){$4n$}
\put(-60,114){$k$}
\put(-402,38){$m$}
\put(-244,38){$4n+1$}
\put(-60,38){$k$}
\caption{A sequence of handleslides and isotopies. Each box \fbox{\hspace{1mm}$\ell$\hspace{1mm}} indicates $\ell $ half-twists; all components without indicated framing carry the framing $-1$.  Applying two isotopies to the top picture yields the second one. From that one the third picture is obtained by sliding each of the two handles containing $2n$ half-twists over the handle to their immediate left (the handles without any twists in them).  The fourth picture is then obtained by a further isotopy and picture 5 is gotten by sliding the handle containing the $4n$ half-twists over the handle to its left. The final picture is gotten by simple additional isotopy. This picture is drawn with $n$ odd in mind; the case of $n$ even works analogously. The components with $m$ and $k$ half-twist are allowed to be the same. }  \label{pic5}
\end{figure}

The surgery description of $Y_m(K)$ is now easily obtained. Namely, in cutting 
$S^3$ open along $D$ (see again figure \ref{pic3}c) one also cuts the $+1$ framed unknot from 
figure \ref{pic3} yielding the $2$-component tangle from figure \ref{pic6}a. To obtain a surgery description of 
$Y_m(K)$, one concatenates $m$ copies of this tangle to obtain an $m$ component link $L=L_1\sqcup ... \sqcup L_m$. The framing $\lambda _i$ of $L_i$ is determined from the equation $1=\lambda _i +\sum_{j\ne i}  \ell k (L_i,L_j)$ (see page 357 of \cite{kauffman2} for an explanation of this formula) and thus turns out to be $\lambda _i = -1$.  Figure \ref{pic6}b illustrates the 
example of $m=4$.

We see that $Y_m(K)$ is obtained by $-1$-framed surgery on an $m$ component \lq\lq necklace\rq\rq, each of whose components is an unknot. To simplify 
this picture, we perform a number of handle slides with the goal of reducing the number of components. The key idea is as depicted in figure \ref{pic5}. Figure \ref{pic5} describes 3 handle slides by which one can  replace 
4 consecutive $-1$-framed components of the \lq\lq necklace\rq\rq by a single component, still with framing $-1$, though this new component is given a single right-handed half-twist. The same holds if one starts with 4 consecutive 
$-1$-framed components that each contain $n$ half-twists (in what we may refer to as the {\em $n$-twisted necklace}): in this case one is left with a single component with framing $-1$ but now with $4n+1$ right-handed half-twists. We shall refer to the sequence of handle slides from figure \ref{pic5}, reducing the number of components from $4$ to $1$, as the {\em reduction procedure}. The reduction procedure can be applied to any ($n$-twisted) necklace of at least 5 components. 

\begin{figure}[htb!] 
\centering
\includegraphics[width=14cm]{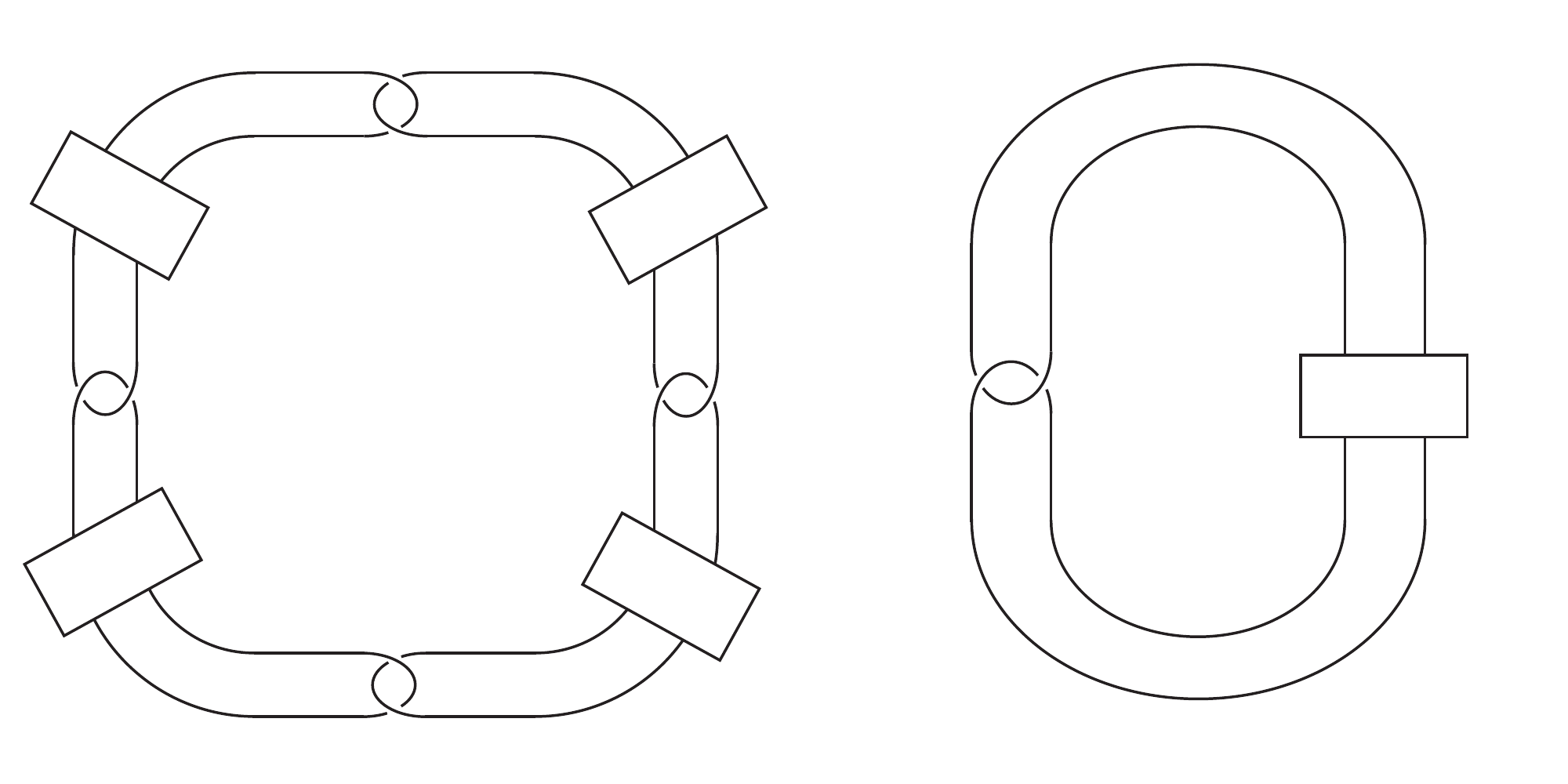}
\put(-305,-5){$(a)$}
\put(-100,-5){$(b)$}
\put(-155,170){$-3$}
\put(-371,144){$\ell$}
\put(-230,142){$\ell$}
\put(-230,46){$\ell$}
\put(-373,52){$\ell$}
\put(-63,95){$4\ell + 1$}
\put(-335,150){$-1$}
\put(-278,150){$-1$}
\put(-278,40){$-1$}
\put(-335,40){$-1$}
\caption{(a) The surgery description of $Y_{4^k}(K)$ obtained by applying the reduction procedure $k-1$ times to the $4^k$ component necklace of $-1$ framed unknots. The number of half-twists $\ell$ in each of the four boxes is given by $\ell = (4^{k-1}-1)/3$. (b) Three additional handle slides turn figure (a) into $-3$-framed surgery on a twist knot. The handle slides performed are almost identical to those from figure \ref{pic5} and are left as an easy exercise.}  \label{pic7}
\end{figure}

We proceed by separately considering $\delta_{2^{2k}}(K)$ and $\delta_{2^{2k+1}}(K)$. To compute the former, 
one uses the reduction procedure described above to reduce the $-1$-framed, $4^k$ component necklace 
describing $Y_{4^k}(K)$, to the $4$ component surgery diagram in figure \ref{pic7}a.  Three further handle slides, akin to those from figure \ref{pic5}, reduce this description to $-3$-framed surgery on a single twist knot, see 
figure \ref{pic7}b.

\begin{figure}[htb!] 
\centering
\includegraphics[width=14cm]{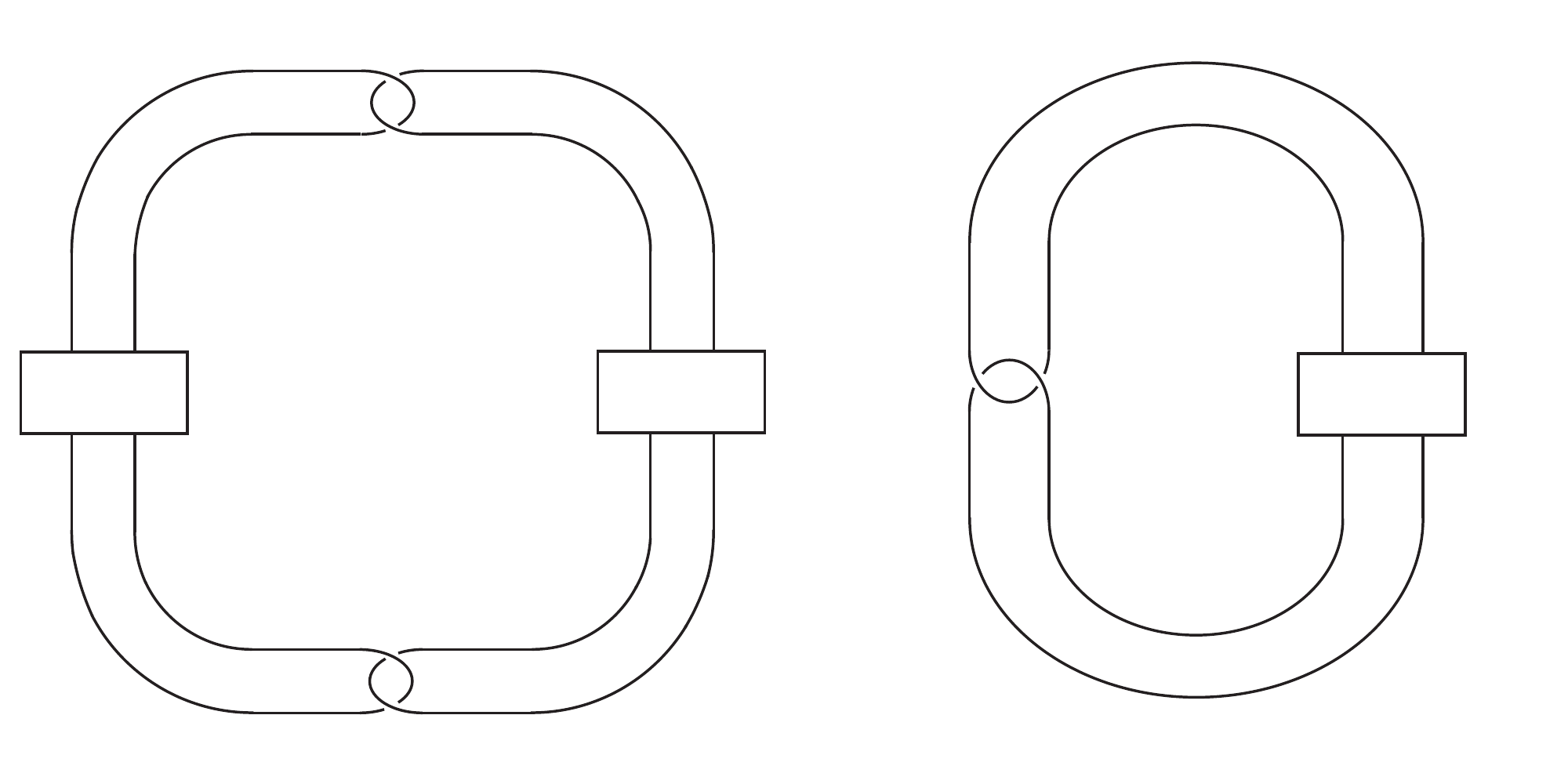}
\put(-305,-5){$(a)$}
\put(-100,-5){$(b)$}
\put(-155,163){$3$}
\put(-229,95){$\ell$}
\put(-375,95){$\ell$}
\put(-54,95){$2\ell$}
\caption{(a) The surgery description of $Y_{2^{2k+1}}(K)$ obtained by applying the reduction procedure $k$ times to the $2^{2k+1}$ component necklace of $-1$ framed unknots. The number of half-twists $\ell$ in each of the two boxes on the left is: $\ell = (4^{k}-1)/3$. (b) An additional handle slides turn figure (a) into $3$-framed surgery on a twist knot.}  \label{pic8}
\end{figure}

\begin{figure}[htb!] 
\centering
\includegraphics[width=6cm]{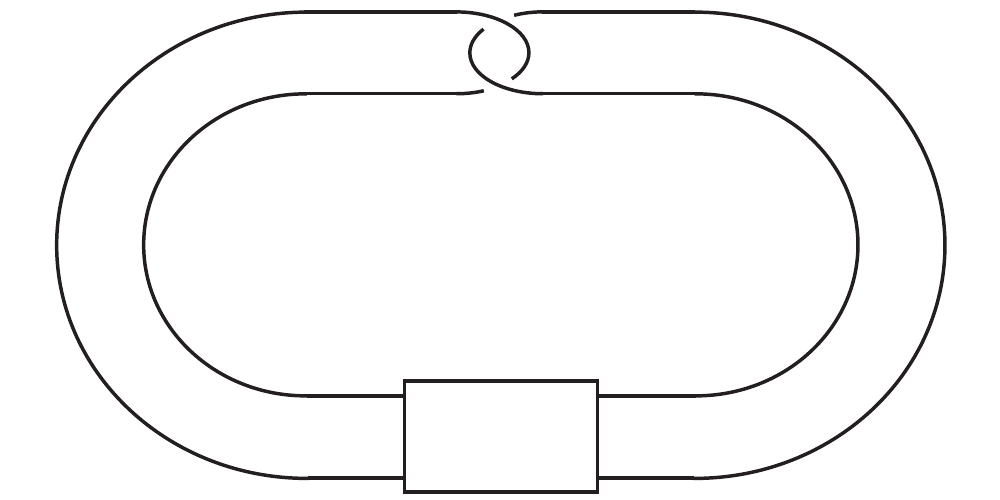}
\put(-88,8){$m$}
\caption{The twist knot $T_m$.}  \label{pic9}
\end{figure}

On the other hand, to find a surgery description of $Y_{2^{2k+1}}(K)$, one repeatedly applies the reduction procedure to the $2^{2k+1}$ component necklace of $-1$-framed unknots to arrive at the $2$-component Dehn 
surgery diagram from figure \ref{pic8}a. A single handle slide (followed by a slam-dunk) gives the description of $Y_{2^{2k+1}}(K)$ from figure \ref{pic8}b.
In summary, we have:
\begin{proposition} \label{surgerieson31}
Let $K$ be the right-handed trefoil and $Y_{2^n}(K)$ be the $2^n$-fold branched cover of $S^3$ with branching set $K$. Let $T_m$ denote the twist knot with the clasp as in figure \ref{pic9}. Then 
$$ Y_{2^n}(K) = \left\{ 
\begin{array}{rl}
-3\mbox{-framed surgery on }T_{(2^n-1)/3}  & \quad ; \quad n=2k, k\ge 1  \cr
& \cr
3\mbox{-framed surgery on }T_{(2^n-2)/3}  & \quad ; \quad n=2k+1, k\ge 0 
\end{array}
\right.
$$ 
\end{proposition}
Using proposition \ref{surgerieson31} and the Heegaard Floer tools for computing correction terms for surgeries on a knot (as described by Ozsv\'ath and Szab\'o in \cite{peter7}), readily lead to a computation of all $\delta_{2^n}$ for the right-handed trefoil.
\begin{theorem}
Let $K$ be the right-handed trefoil knot. Then 
$$ \delta_{2^n}(K) = \left\{
\begin{array}{cl}
3 & \quad ; \quad \mbox{if $n=2k$ and $k\ge 1$.} \cr
& \cr
1 & \quad ; \quad \mbox{if $n=2k+1$ and $k\ge 0$.}
\end{array}
\right.
$$
\end{theorem} 
This theorem is an application of Corollary 4.2 from \cite{peter7}. The only inputs required by that corollary 
are the signatures $\sigma(T_m)$ and Alexander polynomials $\Delta_{T_m}(t)$ of the twists knots. 
These, in turn, are easily determined 
since all $T_m$ are alternating. For any integer $m\ge 1$ one finds: 
$$ \sigma (T_m) = \left\{ 
\begin{array}{cl}
0 & ; \quad m \mbox{ is even} \cr
& \cr
2 & ; \quad m \mbox{ is odd}
\end{array}
\right.
\quad \mbox{ and }\quad 
\Delta_{T_m}(t) = \left\{ 
\begin{array}{cl}
\frac{m}{2}t -( m+1) + \frac{m}{2}t^{-1} & ; \quad m \mbox{ is even} \cr
& \cr
\frac{m+1}{2}t -m + \frac{m+1}{2}t^{-1} & ; \quad m \mbox{ is odd}
\end{array}
\right.
$$ 
Finally, since $|H_1(Y_{2^n}(K);\mathbb{Z})|=3$ for all $n\ge 1$, each $Y_{2^n}(K)$ possesses a unique 
spin-structure which by necessity has to equal $\s_0(K,p^n)$. The corresponding correction term 
$d(Y_{2^n}(K),\s_0(K,p^n))$ is distinguished by the fact that it is the only one, out of the 3 corrections terms of $Y_{2^n}(K)$, which yields an integer when multiplied by 2. This fortuitous coincidence makes the determination 
of $\s_0(K,2^n)$ easy. 

\end{document}